\documentclass{amsart}
\usepackage{citesort}
\usepackage{amssymb,amsmath}
\usepackage{epsfig}
\usepackage{graphicx}
\usepackage{psfrag}
\usepackage{enumerate}

\DeclareMathOperator{\graph}{graph}

\def\R{\mathbb{R}}
\def\N{\mathbb{N}}

\def\theta{\vartheta}
\def\phi{\varphi}
\def\epsilon{\varepsilon}
\def\fracd#1#2{\displaystyle\frac{\displaystyle #1}{\displaystyle #2}}
\long\def\umbruch{{\displaybreak[1]}}

\def\ol#1{\overline{#1}}
\def\dt{{\frac d{dt}}}
\long\def\neueZeile{{\rule{0mm}{1mm}\\[-2ex]\rule{0mm}{1mm}}}

\newtheorem{theorem}{Theorem}[section]
\newtheorem{lemma}[theorem]{Lemma}
\newtheorem{proposition}[theorem]{Proposition}
\newtheorem{corollary}[theorem]{Corollary}

\theoremstyle{definition}

\newtheorem{definition}[theorem]{{Definition}}

\numberwithin{equation}{section}
\newcommand{\abs}[1]{\left\lvert#1\right\rvert}

\begin{document}

\title{Evolution of convex lens-shaped networks under curve shortening
  flow}

\thanks{Several authors are members of and were supported by
  SFB 647/B3 ``{Raum -- Zeit -- Materie:
  Singularity structure, long-time behavior and dynamics of solutions
  of non-linear evolution equations}''}

\def\fuhome{@math.fu-berlin.de}
\address{Oliver C. Schn\"urer, 
  Abderrahim Azouani, 
  Marc Georgi, 
  Juliette Hell,
  Nihar Jangle,
  Tobias Marxen,
  Sandra Ritthaler, 
  Felix Schulze,
  Brian Smith: 
  Freie Universit\"at Berlin, Arnimallee 6, 
  14195 Berlin, Germany}
\author{Oliver C. Schn\"urer}
\email{Oliver.Schnuerer\fuhome}
\author{Abderrahim Azouani}
\email{azouani\fuhome}
\author{Marc Georgi}
\email{georgi\fuhome}
\author{Juliette Hell}
\email{blanca\fuhome\quad(\sc{Juliette Hell})}
\author{Nihar Jangle}
\email{jangle\fuhome}
\author{Amos Koeller}
\address{Amos Koeller: Mathematisches Institut der
Eberhard-Karls-Universit\"at T\"ubingen,
Auf der Morgenstelle 10,
D-72076 T\"ubingen, Germany}
\email{akoeller@everest.mathematik.uni-tuebingen.de}
\author{Tobias Marxen}
\email{marxen\fuhome}
\author{Sandra Ritthaler}
\email{ritthale\fuhome}
\author{Mariel S\'aez}
\address{Mariel S\'aez: Max Planck Institute for Gravitational Physics
(Albert Einstein Institute), Wissenschaftspark Golm, Am M\"uhlenberg 1,
14476 Golm, Germany}
\email{mariel@aei.mpg.de}
\author{Felix Schulze}
\email{Felix.Schulze\fuhome}
\author{Brian Smith}
\email{bsmith\fuhome}

\curraddr{}

\subjclass[2000]{53C44, 35B40}
\date{November 2007.}

\dedicatory{}

\keywords{Mean curvature flow, networks, triple junctions.}

\begin{abstract}
  We consider convex symmetric lens-shaped networks in $\R^2$ that
  evolve under curve shortening flow. We show that the enclosed convex
  domain shrinks to a point in finite time. Furthermore, after
  appropriate rescaling the evolving networks converge to a
  self-similarly shrinking network, which we prove to be unique in an
  appropriate class. We also include a classification result for
  some self-similarly shrinking networks.
\end{abstract}

\maketitle

\markboth{CURVE SHORTENING FLOW OF LENS-SHAPED NETWORKS}{LENS SEMINAR}

\section{Introduction}\label{intro sec}
\noindent
The evolution of networks under curve shortening flow was studied by
Carlo Mantegazza, Matteo Novaga, and Vincenzo Tortorelli
\cite{MantegazzaNetworks}. In their paper they investigated the case
where the underlying graph consists of three lines having one point in
common. In that context they proved short-time existence, long-time
existence and, under extra conditions, convergence to a minimizing
configuration.  However, their methods allowed them to study only tree
networks. In this paper we focus on the evolution of graphs with
closed loops. We consider graphs with special kinds of loops, namely
{\em symmetric lens-shaped networks}.

A convex lens-shaped network $M_0\subset\R^2$ consists 
of two smooth convex arcs and two half-lines
arranged as in the following picture: 
\begin{figure}[htb]
\epsfig{file=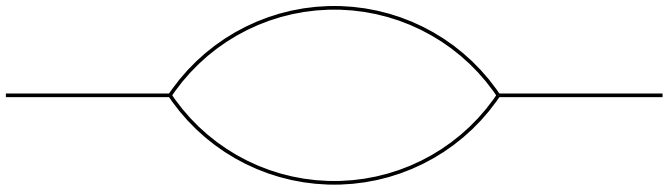,   width=0.7\textwidth}
\caption{Lens-shaped network}
\label{lens pic}
\end{figure}

As in \cite{MantegazzaNetworks}, we assume that the curves meet at the
triple points at an angle of $120^\circ$.  We also impose an extra
symmetry assumption: we consider networks that have reflection
symmetry about the $x^1$-axis.  For precise definitions we refer the
reader to Section \ref{notation sec}.

In order to define the evolution equation of such a network we
consider $G$ to be an abstract graph homeomorphic to $M_0$ consisting
of two bounded intervals and two half-lines. Let $F_0:G\to M_0$ be a
homeomorphism, which restricted to a closed edge is a diffeomorphism.
Now consider a family $(F(\cdot,t))_{t\in I}$ of mappings
$F(\cdot,t):G\to\R^2$, where each $M_t:=F(G,t)$ is a convex
lens-shaped network as described above and $I$ denotes a generic
interval.  Then the family $(M_t)_{0\le t<T}$ is said to evolve under
curve shortening flow if away from the vertices it satisfies
\begin{equation}\label{mcf}
\left(\dt F\right)^\perp=\vec\kappa,
\end{equation}
where $\perp$ denotes the orthogonal projection on the
normal space and $\vec\kappa$ denotes the curvature vector. 
Up to tangential diffeomorphisms, \eqref{mcf} is equivalent
to $\dt F=-\kappa\nu$. Let $\Omega_t$ denote the bounded component of
$\R^2\setminus M_t$. 

We prove that $M_t$ shrinks under curve shortening flow \eqref{mcf} to
a straight line. More precisely,

\begin{theorem}\label{to point thm}
  Let $M_0$ be a convex lens-shaped network. Then there exists a
  family $M_t:=F(G,t)$ of convex lens-shaped networks as above that
  evolves under curve shortening flow \eqref{mcf}. This family is smooth
  for $t>0$ and spatially $C^2$ up to $t=0$.  Moreover,  it exists for
  $t\in[0,T)$, where $T=\frac{3|\Omega_0|}{4\pi}$, and, as $t\nearrow
  T$, the networks $M_t$ converge in Hausdorff distance to the
  $x^1$-axis.
\end{theorem}

This result is analogous to the one proved by Michael Gage and Richard
Hamilton \cite{GageHamilton} for convex curves embedded in $\R^2$.
They showed that such curves shrink to a point in finite time.
Furthermore, after appropriate rescaling, these curves converge to a
circle.  Notice that circles shrinking homothetically at an
appropriate rate are self-similar solutions to \eqref{mcf}. Moreover,
in \cite{AbreschLanger} Uwe Abresch and Joel Langer proved these
circles are the unique embedded self-similarly shrinking solutions to
curve shortening flow in the plane.
 
Similarly, in order to analyze the blow-up profile of shrinking
lenses, we  consider families of networks $(M_t)_{t\in I}$ that
shrink self-similarly. That is, for every $t_1$, $t_2\in I$ the
network $M_{t_1}$ is the image of $M_{t_2}$ under a homothety. The
existence of such networks was studied numerically by
William Mullins \cite{MullinsGrainBdry}. We prove the following
classification theorem:

\begin{theorem}\label{to self sim thm}
  There exists a unique self-similarly shrinking family of symmetric
  lens-shaped networks $(N_t)_{t\in(-\infty,\,0)}$ solving \eqref{mcf}
  such that as $t\nearrow0$ the networks $N_t$ converge to the
  $x^1$-axis in Hausdorff distance.  Moreover, all networks $N_t$ are
  symmetric with respect to the $x^2$-axis.
\end{theorem}

As in \cite{GageHamilton}, this enables us to describe the shape of
the networks in Theorem \ref{to point thm} more closely as $t\nearrow
T$:
\begin{theorem}\label{self sim unique thm}
Let $(M_t)_{t\in[0,T)}$ be a family of networks as in 
Theorem \ref{to point thm}. Then there exists a unique
$x_0$ such that $\{x_0\}=\bigcap\limits_{[0,\,T)}\Omega_t$ and
the rescaled networks converge smoothly to a 
network that contracts self-similarly:
$$\left(2(T-t)\right)^{-1/2}\cdot(M_t-x_0)
\to N_{\left(-\frac12\right)},$$
where $N_t$ is as in Theorem \ref{to self sim thm}.
\end{theorem}

For the sake of completeness, we also include a complete
classification of a certain class of homothetically shrinking
networks. Namely, we consider networks that have the same topology as
the lens-shaped ones. We prove that,  besides the family described in
Theorem \ref{to self sim thm},  there exists precisely one (up to
rotations and reflections)  of such networks.  It is depicted as follows:
\begin{figure}[htb]
\epsfig{file=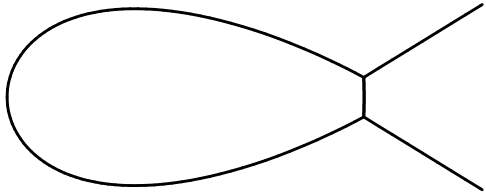,
  width=0.5\textwidth}
\caption{Homothetically shrinking fish-shaped network}
\label{fish pic}
\end{figure}

We organize the paper as follows: In Section \ref{notation sec} we
establish notation, rewrite the equation as a free boundary value
problem and prove short-time existence. Theorem \ref{to point thm} is
proved in Section \ref{long time and point sec}, that is we show
long-time existence up to a time $T$ such that $M_t$ converges to a
line.  In Section \ref{rescaling sec} we carry out a blow-up analysis
that proves Theorem \ref{self sim unique thm}. The existence and
uniqueness of self-similar shrinking lenses is studied in Section
\ref{selfsimilar sec}. Homothetically shrinking fish-shaped networks
are considered in Appendix \ref{fish app}.

Finally, we would like to recall that the result of Michael Gage and
Richard Hamilton was later extended by Matthew Grayson, who proved
that embedded closed curves become convex in finite time
\cite{GraysonRoundPoints}.  Unfortunately, an analogous extension for
our problem is not included in this paper.  However, this might be an
interesting direction to explore in the future.

This paper was written in a research seminar organized by Oliver
Schn\"urer at the Free University Berlin in 2006 and 2007.

\section{Notation}\label{notation sec}

\noindent
A symmetric lens-shaped network $M\subset\R^2$ is determined
by a smooth curve \mbox{$\gamma:[-1,1]\to\R^2$}, which is regular
up to the boundary, and such that the following conditions hold:
\begin{enumerate}[(i)]
\item $\gamma((-1,1))\subset\left\{\left(x^1,x^2\right)\in\R^2:
x^2>0\right\}$,
\item $\gamma(-1),\,\gamma(1)\in\left\{\left(x^1,x^2\right)\in\R^2:
x^2=0\right\}$ and $x^1(\gamma(-1))<x^1(\gamma(1))$,  
\item at the endpoints, the curve and the $x^1$-axis enclose an
angle of $60^\circ$; more precisely,
$$\frac{\gamma'(-1)}{|\gamma'(-1)|}= 
\left(\tfrac{1}{2},\tfrac{\sqrt{3}}{2} \right), \qquad
\frac{\gamma'(1)}{|\gamma'(1)|}=\left(\tfrac{1}{2},
  -\tfrac{\sqrt{3}}{2}\right)\ .$$
\end{enumerate}

Then the complete symmetric lens-shaped network $M$ is given by four
curves: $\gamma_1,\, \gamma_2,\, \gamma_3,\, \gamma_4$, where
$\gamma_1$ is the curve $\gamma$ described above, $\gamma_2$ is the
reflection of $\gamma_1$ about the $x^1$-axis, and $\gamma_3,\,
\gamma_4$ are the half-lines contained in the $x^1$-axis connecting
$\gamma_1(-1)$ and $\gamma_1(1)$ to $-\infty$ and $+\infty$,
respectively. Hence, we can always identify the network $M$ with the
curve $\gamma:=\gamma_1$. Notice that the network is convex (resp.{}
strictly convex) if the curvature of $\gamma$ is non-negative (resp.{}
strictly positive) with an appropriate sign convention. It is easy to
verify that in the interior the family of symmetric lens-shaped
networks $\left(\gamma_t\right)_{0\leq t<T}$ evolves according to
\eqref{mcf} if it satisfies the evolution equation
$$ \langle \dot{\gamma_t}, \nu\rangle = 
\frac{\langle (\gamma_t)_{xx},\nu\rangle}{|(\gamma_t)_x|^2},$$ where
$\nu$ is a choice of unit normal vector to $\gamma_t$.  Here and in
what follows a dot will always denote the derivative with respect to
time and a subscript $x$ the derivative with respect to the spatial
parameter.

For such a family of networks, we define $a(t), b(t) \in \R,\, a(t)<
b(t)$, for $t\in [0,T)$ by $\gamma_t(-1)=\left(a(t),0\right)$ and
$\gamma_t(1)= \left(b(t),0\right)$. Furthermore, let $\Omega_t$ be the
open bounded component of $\R^2\setminus M_t$ and $|\Omega_t|$ its
$2$-dimensional Lebesgue measure. The length of the network is defined
by $L_t:=2\cdot \text{length}(\gamma_t)$.

Suppose that $\gamma_t([-1,1])=\graph u(\cdot,t)|_{[a(t),b(t)]}$ for a
function $u:\R\times[0,T')\rightarrow [0,\infty)$.  Then we call the
corresponding network graphical.  Notice that for times when a network
is convex, it is also graphical. Moreover, the function $u$ is concave
in $[a(t),b(t)]$ and the network remains graphical for some time.  In
this situation the problem can be reformulated as follows:

A continuous function $u:\R\times[0,T)\to\R$ determines
a solution to \eqref{mcf},
if there exist $a:[0,T)\to\R$ and $b:[0,T)\to\R$ with
$a(t)<b(t)$ for all $t\in[0,T)$ and
\begin{equation}\label{free bvp}
\begin{cases}
\dot u=\frac{\displaystyle u_{xx}}{\displaystyle 1+u_x^2}&
\text{in the interior of }
\bigcup\limits_{t\in[0,T)}(a(t),b(t))\times\{t\},\\
u(x,t)=0&\text{for }x\in(-\infty,a(t)]\cup[b(t),\infty)
\text{ and }t\in[0,T),\\
u_x(a(\cdot),\cdot)=\tan\frac\pi3=\sqrt3&\text{in }[0,T),\\
u_x(b(\cdot),\cdot)=-\tan\frac\pi3=-\sqrt3&\text{in }[0,T),\\
u(\cdot,0)=u_0&\text{in }\R,
\end{cases}
\end{equation}
where $u_0:\R\to\R$ is such that 
\begin{equation}
\begin{cases}
u_0=0&\text{in }(-\infty,a(0)]\cup[b(0),\infty),\\
(u_0)_x(a(0))=\tan\frac\pi3=\sqrt3,&\\
(u_0)_x(b(0))=-\tan\frac\pi3=-\sqrt3.&\\
\end{cases}
\end{equation}
Here expressions like $u_x(a)$ have to be considered as the limits
from above, whereas expressions like $u_x(b)$ are limits from below.
The solutions considered are in $C^{2;1}$ and, hence,  
we require that $u_0\in C^2([a(0),b(0)])$.  Furthermore, we choose 
$u_0$ such that $\graph u_0\cup\graph(-u_0)=M_0$.

Under this formulation it is easy to see that the maximum principle
implies $u(x,t)>0$ for $a(t)<x<b(t)$ and that any solution to
\eqref{mcf} is unique in the class of lens-shaped networks.

Short-time existence for \eqref{mcf}, and thus also for \eqref{free
  bvp}, follows from \cite[Theorem 3.1]{MantegazzaNetworks} by
introducing artificial boundary points on the half-lines. Furthermore,
estimates in \cite{MantegazzaNetworks} imply that solutions are in
$C^{2;1}$.  Summarizing we have

\begin{theorem}[\cite{MantegazzaNetworks}]\label{C21 thm}
For any smooth symmetric lens-shaped network, \eqref{mcf} has a 
solution of class $C^{2;1}$ up to the free boundary 
on a short time interval.\end{theorem}

We finish this section by writing the relevant geometric quantities in
terms on the function $u$. We also recall the definition of the
support function, that will be  used in Section \ref{selfsimilar sec}.

An easy computation shows that the upwards pointing unit normal of 
$\gamma_t=\graph u(\cdot,t)|_{[a(t),b(t)]}$ is  
$$\nu=\frac{(-u_x,1)}{\sqrt{1+u_x^2}}.$$  
The curvature is  given by  
$$\kappa=\frac{-u_{xx}}{(1+u_x^2)^{3/2}}$$ 
and as usual the curvature vector is given by $\vec\kappa=-\kappa\nu$.
We define $\nu$ along the $\graph(-u(\cdot,t)|_{[a(t),b(t)]})$ as the
reflection of the normal vector $\nu$ from above. On $\{x^2<0\}$ we
define $\kappa$ such that $\kappa\ge0$ for a convex lens. 

We shortly recall that for a curve $\alpha:I\to\R^2$, the
support function $S:\nu(\alpha(I))\to\R$ is defined by
$S(p):=\langle\alpha(\nu^{-1}(p)),p\rangle$. More details and
properties of the support function can be found in Section
\ref{selfsimilar sec}.

\section{Long-Time Existence and Convergence to a Point}
\label{long time and point sec}

In this section we prove long-time existence and convergence to a
point. We do this via a series of lemmata. The main ingredients are
gradient estimates, a formula for the rate of decrease of the volume
and curvature bounds. These last estimates are proven using a modified
intrinsic-extrinsic distance ratio, that was introduced by Gerhard
Huisken in \cite{HuiskenAsian}.
 
Note that for the first two lemmata we do not have to assume that our
networks are convex.
\begin{lemma}[Solutions stay graphical]
Let $M_0$ be a graphical symmetric lens-shaped network. If
\eqref{mcf} has a solution on a time interval $[0,T)$ then
$M_t$ is graphical for all $t\in[0,T)$. In fact, we have
for $u$ as in \eqref{free bvp} the estimate
\begin{equation}\label{C1 est}
|u_x|\le\max\left\{\sqrt3,\,\sup\limits_{x\in(a(0),b(0))}
|u_x(x,0)|\right\}.
\end{equation}
\end{lemma}
\begin{proof}
  For a smooth network being graphical means that $|u_x|<\infty$.
  Notice that for a smooth evolution of networks this condition is
  preserved for short times. Moreover, in order to prove the lemma it
  suffices to consider \eqref{free bvp} instead of \eqref{mcf} and to
  show that $t\mapsto\max\limits_{x\in[a(t),b(t)]}|u_x(x,t)|$ is
  non-increasing.

We employ the maximum principle on
the domains where $u>0$, 
that is for each $t$ we take $x\in (a(t),b(t))$.
 On $\{(a(t),t):t\in[0,T)\} \cup
\{(b(t),t):t\in[0,T)\}$, we have $|u_x|=\sqrt3$. 
It is easy to compute that the
evolution equation for $u_x$ is given by
$$\dt(u_x)=\frac1{1+u_x^2}(u_x)_{xx}
-\frac{2u_xu_{xx}}{\left(1+u_x^2\right)^2}(u_x)_x.$$ Hence, according
to the maximum principle, $|u_x|$ cannot attain a new maximum in the
interior of $\bigcup\limits_{t\in(0,T)}(a(t),b(t))\times\{t\}$.
\end{proof}

\begin{lemma}[Evolution of the enclosed volume]
\label{graphical volume lem}\rule{1ex}{0em}
Let $(M_t)_{t\in[0,T)}$ be a family of symmetric networks 
solving \eqref{mcf}. 
Then the enclosed volume fulfills 
$$\dt|\Omega_t|=-\frac{4\pi}3.$$
\end{lemma}
\begin{proof}
  
  We give a proof only for graphical networks. In particular, our
  proof shows that there is no contribution to the change of enclosed
  area from the triple points. Near the triple points, a smooth
  symmetric network is always graphical. Thus, the general result can
  be obtained similarly as for closed curves \cite{GageHamilton}.
\par
For a graphical network, we have
$$|\Omega_t|=2\int\limits_{a(t)}^{b(t)}u(x,t)\,dx.$$
Differentiating yields
\begin{align*}
  \dt|\Omega_t|=&2\underbrace{u(b(t),t)}_{=0}\dot b(t)
  -2\underbrace{u(a(t),t)}_{=0}\dot a(t)
  +2\int\limits_{a(t)}^{b(t)}\dot u(x,t)\,dx\umbruch\\
  =&2\int\limits_{a(t)}^{b(t)}\frac{u_{xx}}{1+u_x^2}(x,t)\,dx
  =2\arctan(\underbrace{u_x(b(t))}_{=-\sqrt3})
  -2\arctan(\underbrace{u_x(a(t))}_{=\sqrt3}) =-\frac{4\pi}3.
\end{align*}  
The lemma follows.
\end{proof}

Note that this lemma also implies the following: If $\Omega_t$ shrinks
to a point at time $T$ or, more generally,
$\bigcap\limits_{t\in[0,T)}\Omega_t$ is 
a set of measure zero, then 
$T=\frac{3\cdot|\Omega_0|}{4\pi}$.

\begin{lemma}[Convexity is preserved]\label{pres conv lem}
A convex network evolving (smoothly) 
according to \eqref{free bvp} stays convex.
More precisely, a non-negative lower bound on the
curvatures (in the non-flat part of the network) 
is preserved during the evolution.
\end{lemma}
\begin{proof}
  Let us first assume that the curvature of $\graph u(\cdot,0)$ is
  uniformly bounded below in $(a(0),b(0))$ by $\lambda>0$ (that is,
  the network is strictly convex).
  
A simple computation in \cite{GageHamilton} shows that the
evolution equation for the curvature $\kappa$ is
$$\dt\kappa=\Delta\kappa+\kappa^3,$$
where $\Delta\kappa$ denotes the Laplace-Beltrami operator applied to
$\kappa$. 

Let
$$k(t):=\inf\limits_{\genfrac{}{}{0pt}{}{\tau\in[0,t]}{x\in[a(t),b(t)]}}
\kappa(x,\tau).$$ Note that $k(t)$ is monotonically decreasing,
Lipschitz continuous and hence differentiable for almost every $t$.
Moreover, the evolution equation for $\kappa$ and the maximum
principle imply that at times that $k(t)$ is decreasing and
differentiable, the minimum has to be attained at the boundary and the
derivative of $k(t)$ equals $\dt\kappa(a(t),t)$ or
$\dt\kappa(b(t),t)$.

Fix $0<\epsilon<\lambda$ arbitrarily and suppose that there is a $t_1$
such that $k(t_1)=\lambda-\epsilon$. Assume that $t_1$ is minimal with
that property. Hence there exists a time $t_0$ with $0<t_0<t_1$ such
that $\dt k(t_0)<0$. Without loss of generality we may assume that
$\dt k(t_0)=\dt\kappa(a(t_0),t_0)$ and get
\begin{align}\label{def t0 conv pres}
\left.\dt\kappa(a(t),t)\right|_{t=t_0}<&\ 0,\umbruch\\
\kappa(a(t_0),t_0)=&\ \lambda-\tilde\epsilon,\quad
0<\tilde\epsilon<\epsilon,\nonumber\umbruch
\intertext{and}\umbruch
\kappa(a(t_0),t_0)\le&
\inf\limits_{\genfrac{}{}{0pt}{}{t\in[0,t_0]}{a(t)\le x\le b(t)}}
\kappa(x,t).\nonumber
\end{align}

By rotating our network $90^\circ$ we can describe it near the triple
point $(a(t),0)$ as $\graph\tilde u|_{[0,\delta)\times[0,\tau)}$ for
some $\delta>0$ and $\tau>t_0$. Thus $\tilde u$ solves
\begin{align}\label{rotated free bvp}
\begin{cases}
  \dot{\tilde u}=\frac{\tilde u_{xx}}{1+\tilde u_x^2}&\text{in
  }[0,\delta)  \times[0,\tau),\\
  \tilde u_x(0,t)=-\frac1{\sqrt3}&\text{in }[0,\tau).
\end{cases}
\end{align}

Now we construct a barrier for $\tilde{u}$ using the downwards
translating grim reaper, given by
$$\graph\left(\left(-\frac\pi2,\frac\pi2\right)\ni x\mapsto
  \log\cos x\right).$$ First notice that by rescaling and translating
the grim reaper it is possible to get a $\zeta>0$ and a function
$\Gamma:[0,\zeta)\times\R\to\R$ that satisfies
$$\begin{cases}
  \dot\Gamma=\frac{\Gamma_{xx}}{1+\Gamma_x^2}&\text{in
  }[0,\zeta)\times\R,\\
  \Gamma_x(0,t)=-\frac1{\sqrt3}&\text{for all }t\in\R,\\
  \frac{-\Gamma_{xx}}{\left(1+\Gamma_x^2\right)^{3/2}}(0,t)
  =\lambda-\tilde\epsilon
  &\text{for all }t\in\R,\\
  \Gamma(0,t_0)=\tilde u(0,t_0),&\\
  \Gamma(x,t)\to-\infty&\text{as }x\nearrow\zeta\text{ for every
  }t\in\R.
\end{cases}$$

\begin{figure}[htb]
\psfrag{u}{$\graph u$}
\psfrag{grim reaper}{grim reaper}
\epsfig{file=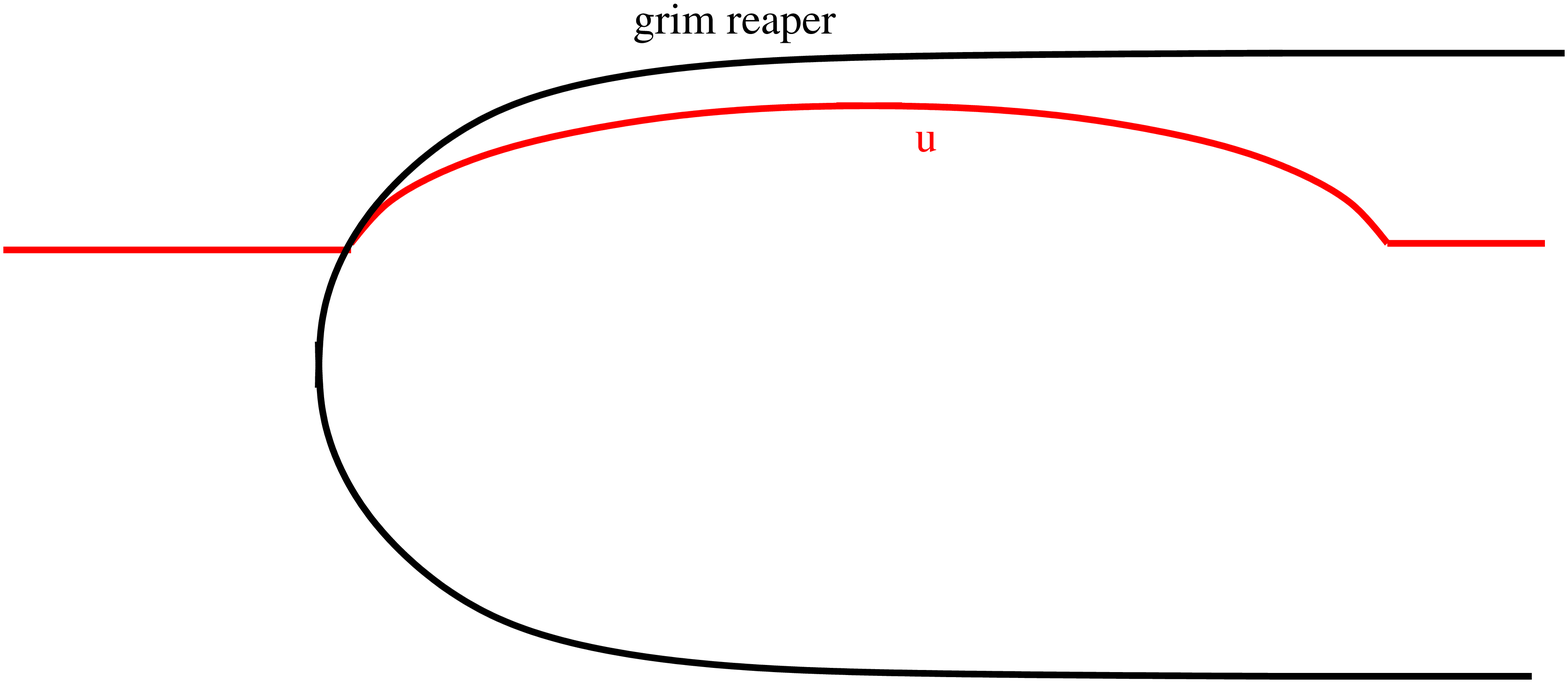, width=0.5\textwidth}
\end{figure}

As the curvature of $\graph\Gamma(\cdot,t_0)$ is bounded above by
$\lambda-\tilde\epsilon$ and the curvature of $M_{t_0}$ is bounded
below by $\lambda-\tilde\epsilon$ (with equality at the triple point),
we see that $\delta<\zeta$ and 
\begin{equation}\label{barrera 1}
  \Gamma(x,t_0)\ge\tilde u(x,t_0) \hbox{ for }
  x\in[0,\delta),\end{equation}  with strict inequality  near
$x=\delta$. From the 
maximum principle follows that \begin{equation}\label{barrera 1.5}
  \Gamma(x,t)\ge\tilde u(x,t) \hbox{ for all }
  (x,t)\in[0,\delta)\times[t_0,\tau).\end{equation}  Moreover, our
construction implies that 
the curvatures of 
$\graph\tilde u(\cdot,t_0)$ and of $\graph\Gamma(\cdot,t_0)$ 
coincide for $x=0$. Thus the evolution equation and
the boundary condition imply that 
\begin{equation}\label{barrera 2}\dot\Gamma(0,t_0)=
\dot{\tilde u}(0,t_0)=\frac{\tilde u_{xx}}{\left(1+\tilde
    u_x^2\right)^{3/2}} \sqrt{1+\tilde u_x^2}(0,t_0)
=-\frac2{\sqrt3}\kappa(0,t_0).\end{equation}

Combining \eqref{barrera 1}, \eqref{barrera 1.5} and \eqref{barrera 2}
we conclude that $\ddot\Gamma(0,t_0)\ge\ddot{\tilde u}(0,t_0)$ and
thus
$$\dot\kappa(0,t_0)=-\frac{\sqrt3}2
\ddot{\tilde u}(0,t_0) \ge-\frac{\sqrt3}2 \ddot\Gamma(0,t_0).$$ Since
$\graph\Gamma(\cdot,t)$ is a translating solution we have
$\ddot\Gamma=0$. This fact combined with the equation above implies
$\dot\kappa(0,t_0)\ge0$, which contradicts \eqref{def t0 conv pres}
and finishes the proof in the strictly convex case.

Now the general case follows by approximation. Indeed every convex
initial condition may be approximated in $C^2$ by strictly convex
networks.  Standard results (see also \cite{MantegazzaNetworks} and
Theorem \ref{C21 thm}) imply that solutions for this Neumann problem
are in $C^{2;1}$ and depend continuously on the initial data.
Moreover, by our previous proof we have that the approximations stay
strictly convex and $C^2$-close to our network for some positive time
interval.  The lemma follows by taking limits.
\end{proof}

In what follows we consider the reduced network $G'$ consisting of the
network $G$ without the two open half-lines. We often identify $G'$
with $F(G',t)$. The coming lemmata prove that along the flow
$(F(\cdot,t))_{t \in [0,T)}$ (always restricted to $G'$) the ratio of
extrinsic and (modified) intrinsic distance defined in
\cite{HuiskenAsian} is bounded from below.

Recall that the extrinsic distance $d_{\text{ex}}$ is given by
\[ d_{\text{ex}}: G' \times G' \times [0,T) \to \R, \,\,\,
d_{\text{ex}}(p,q,t) = \|F(p,t) - F(q,t)\|_2, \] where $\|\cdot\|_2$
denotes the Euclidean norm on $\R^2$.  The intrinsic
distance $d_{\text{in}}$ can be computed as
\[ d_{\text{in}}: G' \times G' \times [0,T) \to \R, \,\,\,
d_{\text{in}}(p,q,t) = \int\limits_p^q ds_t, \] where $s_t$ is the
arc-length parameter at time $t$. Note that there are two locally
minimizing paths from $p$ to $q$ in $F(G',t)$ of lengths smaller than
the total length of $F(G',t)$. Let $d_{\text{in}}$ be the length of
the shorter path $\gamma_{\text{in}}$, which is less or equal than
$L/2$.  Then the modified intrinsic distance is given by $\psi =
\frac{L}{\pi} \sin\left(\frac{\pi d_{\text{in}}}{L}\right)$. Note that
$\psi$ is smooth if $d_{\text{in}}>0$.
\begin{proposition}
  Consider a family of convex symmetric networks solving \eqref{mcf}
  (or equivalently \eqref{free bvp}). Suppose that at time $t_0 \in
  [0,T)$ the quotient $\frac{d_{\text{ex}}}{\psi}$ attains a spatial
  minimum at $(p,q)$ such that $p\neq q$ and neither $p$ nor $q$ is a
  triple point. Then
  \[ \left.\frac{d}{dt}
    \frac{d_{\text{ex}}}{\psi}(p,q,t)\right|_{t=t_0} \geq 0. \]
\end{proposition}
\begin{proof}
  We reduce the calculation to the one done in \cite{HuiskenAsian} in
  the regular case.  In order to do this we show that the terms
  obtained when computing $ \frac{d}{dt}\frac{d_{\text{ex}}}{\psi}$
  are equal to the ones in \cite{HuiskenAsian} plus extra positive
  contributions coming from the triple points.  Let us first compute
  $\frac{d}{dt} L = \dot{L}$, $\frac{d}{dt} d_{\text{in}} =
  \dot{d_{\text{in}}}$ and $\frac{d}{dt}\psi = \dot{\psi}$.  Since
  \[ L(t) = 2 \int\limits_{a(t)}^{b(t)}ds_t = 2
  \int\limits_{a(t)}^{b(t)}\sqrt{1 + u_x^2}\,dx,\] 
  we have
\begin{align*}
  \dt L(t) =& 2 \int\limits_{a(t)}^{b(t)}\frac{u_x \dot{u}_x}{\sqrt{1 +
      u_x^2}}\,dx + 2\sqrt{1 + u_x^2(b(t),t)}\cdot \dot{b}(t) -
  2\sqrt{1 + u_x^2(a(t),t)}\cdot \dot{a}(t)\umbruch\\ 
  =& 2 \int\limits_{a(t)}^{b(t)} \frac{u_x}{\sqrt{1 + u_x^2}}
  \left(\frac{u_{xx}}{1 + u_x^2}\right)_x \,dx + 4\dot{b}(t) -
  4\dot{a}(t).
\end{align*}
By differentiating $u(a(t),t)=0$ and $u(b(t),t)=0$ with respect to $t$
and using \eqref{free bvp} we obtain $u_{xx}(a(t),t) =
-4\sqrt{3}\dot{a}(t)$ and $u_{xx}(b(t),t) = 4\sqrt{3}\dot{b}(t)$.  Now
integrating by parts it follows that
\begin{align*}
  \dt L(t) =& -2 \int\limits_{a(t)}^{b(t)} \left(\frac{u_x}{\sqrt{1 +
        u_x^2}}\right)_x \frac{u_{xx}}{1 + u_x^2} \,dx + \dot{b}(t) -
  \dot{a}(t)\umbruch\\
  =& -2 \int\limits_{a(t)}^{b(t)} \kappa^2 \sqrt{1 + u_x^2} \,dx +
  \dot{b}(t) - \dot{a}(t)\umbruch\\
  =& -2 \int\limits_{(a(t),0)}^{(b(t),0)} \kappa^2 \,ds_t +
  \left[\dot{b}(t) - \dot{a}(t)\right].
\end{align*}
Here we used $\kappa^2 = \frac{u_{xx}^2}{(1 + u_x^2)^3}$. The first
term appears in the computations in \cite{HuiskenAsian} while the
second is a consequence of the triple points.

 If neither $p$
nor $q$ is a triple point a similar calculation implies
\[ \dt d_{\text{in}} = \begin{cases} -\int\limits_p^q \kappa^2 ds_t -
  \dot{a} & \text{if }(a(t),0)\in\gamma_{\text{in}},\\
  -\int\limits_p^q \kappa^2 ds_t + \dot{b} &
  \text{if }(b(t),0)\in\gamma_{\text{in}},\\
  -\int\limits_p^q \kappa^2 ds_t & \text{otherwise.}
\end{cases} \]
The chain rule and the computations  above yield
\begin{align*}
\dt \frac{d_{\text{ex}}}{\psi} =& \frac{\dot{d_{\text{ex}}}}{\psi} -
\frac{d_{\text{ex}}}{\psi^2}\dot{\psi}\umbruch\\ 
=& \frac{\dot{d_{\text{ex}}}}{\psi} - \frac{d_{\text{ex}}}{\psi^2}
\left(\frac{\dot{L}}{\pi}\sin\left(\frac{\pi d_{\text{in}}}{L}\right)
  + \dot{d_{\text{in}}} \cos\left(\frac{\pi d_{\text{in}}}{L}\right) -
  \frac{\dot{L}d_{\text{in}}}{L} \cos\left(\frac{\pi
      d_{\text{in}}}{L}\right)\right)\umbruch\\ 
=& \frac{\dot{d_{\text{ex}}}}{\psi} -
\frac{d_{\text{ex}}}{\psi^2}\left(-\frac{1}{\pi}\int\limits_{G'}\kappa^2
  ds_t\cdot \sin\left(\frac{\pi d_{in}}{L}\right) +
  \frac{d_{in}}{L}\int\limits_{G'}\kappa^2 ds_t\cdot \cos\left(\frac{\pi
      d_{\text{in}}}{L}\right)\right.\\ 
&\qquad\qquad\qquad\left. - \int\limits_p^q\kappa^2 ds_t\cdot
  \cos\left(\frac{\pi d_{\text{in}}}{L}\right)\right)+ B + C,
\end{align*}
where
\[ B = -\frac{d_{\text{ex}}}{\psi^2} \left(\frac{1}{\pi}\left(\dot{b}
    - \dot{a}\right)\sin\left(\frac{\pi d_{in}}{L}\right) -
  \frac{d_{in}}{L}\left(\dot{b} - \dot{a}\right)\cos\left(\frac{\pi
      d_{\text{in}}}{L}\right)\right) \] and  
\[ C = \begin{cases} \frac{d_{\text{ex}}}{\psi^2}
  \dot{a}\cos\left(\frac{\pi d_{\text{in}}}{L}\right) & \text{if
  }(a(t),0)\in\gamma_{\text{in}},\\
-\frac{d_{\text{ex}}}{\psi^2} \dot{b}\cos\left(\frac{\pi
    d_{\text{in}}}{L}\right) & \text{if
  }(b(t),0)\in\gamma_{\text{in}},\\
0 & \text{otherwise.}
\end{cases}\] Notice that besides $B$ and $C$ all the terms in the
above computation appear as well in the regular case. Therefore, in
order to conclude the result it is enough to show that $B$, $C\ge0$.
Because the lens is convex we have $\dot{a} \geq 0$ and $\dot{b} \leq
0$. Furthermore $0 \leq d_{\text{in}} \leq L/2$, which implies $C \geq
0$.  Similarly, recalling that for $y\ge0$ it holds that $\tan y \geq
y$, we get
\[ B\cdot \frac{\psi^2}{d_{ex}} = \frac1\pi
\left(\dot{a}-\dot{b}\right)\left(\sin\left(\frac{d_{\text{in}}\pi}
    {L}\right) -\frac{d_{\text{in}}\pi}{L} \cos
  \left(\frac{d_{\text{in}}\pi}{L}\right)\right) \geq 0. \]  
\end{proof}
\noindent This result together with the following lemma gives us the
desired bound from below on the ratio of extrinsic and (modified)
intrinsic distance.
\begin{lemma} \label{lem:dextintr.bound} Consider a family of convex
  symmetric networks solving \eqref{mcf} (or equivalently \eqref{free
    bvp}). If at a time $t \in [0,T)$ the quotient
  $d_{\text{ex}}/\psi$ attains a spatial minimum smaller than
  $\sqrt{3}/2$ at $(p,q)$, then neither $p$ nor $q$ can be a triple
  point. Thus we have
$$ \inf\limits_{p\neq q} \frac{d_{\text{ex}}(p,q,t)}{\psi(p,q,t)}\geq 
\min\left\{\inf\limits_{p\neq q}
  \frac{d_{\text{ex}}(p,q,0)}{\psi(p,q,0)},
  \frac{\sqrt{3}}{2}\right\}$$ for all $t \in [0,T)$.
\end{lemma}
\begin{proof}
  We fix $t$. Near the triple points the network is approximately
  straight. Hence, geometric considerations show that
$$\liminf\limits_{p,q}\frac{d_{\text{ex}}(p,q,t)}{\psi(p,q,t)}
=\frac{\sqrt{3}}{2}$$ if $p$ and $q$ approach a triple-point from two
different sides. If they approach another point or the triple point
from the same side, the limit is $1$.\par
If $p$ and $q$ are both triple points it is easy to see that we can
decrease $d_{\text{ex}}/\psi$ by moving $p$ and $q$ clockwise (or
counterclockwise) along the lens maintaining the intrinsic distance
$d_{\text{in}}(p,q,t)$ (and thus maintaining the value of
$\psi(p,q,t)$). So we restrict ourselves to the case that $p$ is the
left triple-point and $q$ lies on the upper part of the lens. The idea
here is again to move $p$ and $q$ either clockwise or
counter-clockwise along the lens maintaining $d_{\text{in}}(p,q,t)$
and decreasing $d_{\text{ex}}(p,q,t)$.\par
Let $\vec{n}^+$ and $\vec{n}^-$ denote the two unit tangent vectors at
the left triple point $p$, pointing to the right and into the upper
half-plane and into the lower half-plane, respectively. Let
$\vec{\tau}$ denote the unit tangent vector at $q$ pointing towards
$p$. Since the lens is convex, we see that if $\vec{\tau}$ is not
parallel to $\vec{n}^-$, we decrease $d_{\text{ex}}(p,q,t)$ by moving
both points counter-clockwise as $\langle\vec n^-,\vec\tau\rangle<0$.
If $\vec{\tau}$ is parallel to $\vec{n}^-$ we can decrease
$d_{\text{ex}}(p,q,t)$ by moving both points clockwise.
\end{proof}

Since $d_\text{in}\leq L/2$ we have $\psi \geq
\frac{2}{\pi}d_\text{in}$. This implies

\begin{corollary}[Bounds on intrinsic and extrinsic distances]
\label{din dex cor} 
 Consider a family of convex
  symmetric networks solving \eqref{mcf}. We have
$$ \frac{d_{\text{ex}}(p,q,t)}{d_\text{in}(p,q,t)}\geq 
\frac{2}{\pi}\min\left\{\inf\limits_{r\neq s}
  \frac{d_{\text{ex}}(r,s,0)}{\psi(r,s,0)},
  \frac{\sqrt{3}}{2}\right\}$$ 
for all $p,q \in G^\prime, p \neq q$ and $t \in [0,T)$.
\end{corollary}

We finish this section by showing:

\begin{proposition}[Curvature Bounds]\label{curvature bounds prop}
  Let $(M_t)_{t\in[0,T)}$ be a family of convex symmetric networks
  solving \eqref{mcf} (or, equivalently, \eqref{free bvp}). Then the
  curvature of $M_t$ is bounded everywhere in terms of a positive
  lower bound for $T-t$, $\sup\kappa(\cdot,0)$, and a bound on the
  quotient of intrinsic and extrinsic distances as obtained in Corollary
  \ref{din dex cor}.
\end{proposition}

Theorem \ref{to point thm} is a 
 direct consequence of the curvature bounds. That  is,  we may 
choose $T$ in \eqref{free bvp} such that $|\Omega_t|\searrow0$
as $t\nearrow T$, i.\,e.\ a solution exists until the enclosed
volume goes to zero. 
\par
We would like to remark that 
from the proof of  Proposition \ref{curvature bounds prop} 
we obtain,  as long as $|\Omega_t|$
is uniformly bounded below, 
uniform
curvature estimates for the sequence of rescaled solutions  
considered in Section \ref{rescaling sec}. 

\begin{proof}[Proof of Proposition \ref{curvature bounds prop}]
  We organize the coming proof as follows. We start by controlling the
  geometry of the network, using the previous bounds on
  $\frac{d_{\text{ex}}}{d_{\text{in}}}$. More precisely, we conclude
  that $\sup|u(\cdot,t)|$, $\sqrt{|\Omega_t|}$ and $|a(t)-b(t)|$ are
  comparable. This allows us to use interior estimates to bound the
  curvature at a given distance away from the triple points. We
  conclude by controlling the curvature up to the triple points via the
  maximum principle.\par
  \textbf{Geometric control:} Assume that $0<\frac1c\le T-t\le c$. From
  Lemma \ref{graphical volume lem} we get similar bounds for
  $|\Omega_t|$. As usual, in view of the short-time existence, we may
  assume that we have a smooth solution up to a certain time, up to
  which we want to prove a priori estimates.\par
  We claim that there exists a constant $c$ depending only on
  $\frac{d_{\text{in}}}{d_{\text{ex}}}$, such that
\begin{align*}
\frac1c\cdot|a(t)-b(t)|\le&\sup\limits_{x\in(a(t),b(t))}|u(x,t)|
\le c\cdot|a(t)-b(t)|
\intertext{and}
\frac1c\cdot|\Omega_t|\le&|a(t)-b(t)|^2\le c\cdot|\Omega_t|.
\end{align*}
These can be seen as follows: Convexity implies the gradient estimate
$-\sqrt3\le u_x(x,t)\le\sqrt3$ for all $x$, $t$. From this follows the
upper bound in the first inequality. We prove the lower bound by
contradiction. If there is a sequence $t_i$ such that
$\frac{\sup|u(\cdot,t_i)|} {a(t_i)-b(t_i)}\to0$, the lens-shaped
domain would be very thin contradicting the upper bound on
$\frac{d_{\text{in}}}{d_{\text{ex}}}$.\par Using convexity again, we
see that $|\Omega_t|$ is comparable to
$$|a(t)-b(t)|\cdot\sup\limits_{x\in(a(t),b(t))}|u(x,t)|,$$ which
implies the second inequality. 
\par \textbf{Curvature bounds:} The short-time existence proof implies
the claimed estimates as long as $t$ is small in terms of
$\sup|u(\cdot,0)|$, $\frac{d_{\text{in}}}{d_{\text{ex}}}$ and
$\sup\kappa(\cdot,0)$.\par
For $t\ge\frac1c>0$, we will use the following interior estimate,
obtained from
\cite{JCcontMCF,ColdingMinicozzi,EckerHuiskenInvent,EckerHuiskenAnn}:
Let $u:(-2r,2r)\times(0,4r^2)\to\R$ for some $r>0$ be a solution to
$\dot u=\frac{u_{xx}}{1+u_x^2}$. Then there exist constants $c_k$,
that depend on $k$, $r$ and
$\Vert u\Vert_{C^0((-2r,2r)\times(0,4r^2))}$, such that
\begin{equation}\label{int est ii}
\sup\limits_{(x,t)\in(-r,r)\times\left(3r^2,4r^2\right)}
\left|\frac{\partial^ku}{\partial x^k}\right|(x,t)\le c_k.
\end{equation}
We would like to remark that in our situation even the weaker
estimates in \cite{EckerHuiskenInvent} would suffice.
\par 
Fix $0<\epsilon<\frac14$ and define
$$M:=\inf\limits_t\sup\limits_x|u(x,t)|.$$
Assume that we are in the set where $u\ge\epsilon\cdot M$. Since
$|u_x|\le\sqrt3$, we can find $\delta>0$ such that this set is
contained in $[a(t)+\delta,b(t)-\delta]\times\R$.  Note that $\dot
a\ge0$ and $\dot b\le0$.  Thus we can apply the interior estimates
\eqref{int est ii} and get the claimed a priori estimates in the set
where $u\ge 2\epsilon\cdot M$.
\par Suppose now that we are in the domain where $u\le
(1-\epsilon)\cdot M$. In this set, convexity implies that $|u_x|$ is
bounded below by $\frac{\epsilon\cdot M}{a(t)-b(t)}$. Thus, rotating
the coordinate system by $\frac\pi2$, $\graph u|_{(a(t),b(t))}$ is
represented in each connected component of this domain as a graph with
bounded gradient.
\par Consider graphs as in \eqref{rotated free bvp}.
Due to the curvature estimates away from the triple point obtained
above, we may assume that we have $C^2$ bounds if $x$ is bounded away
from zero, i.\,e.\ especially for some fixed positive $x$.  Consider
$V:=-\dot u=-\frac{u_{xx}}{1+u_x^2}$ corresponding to $\kappa v$ in
standard covariant notation.  We have the evolution equation
$$\dt V
= \frac1{1+u_x^2}V_{xx}
-\frac{2u_xu_{xx}}{\left(1+u_x^2\right)^2}V_x.$$ Differentiating
the boundary condition yields $V_x=0$ at the triple point. Now the
claim follows from the maximum principle.
\end{proof}

\section{Proof of Theorem \ref{to self sim thm}}\label{rescaling sec}
\noindent
In this section we combine the bounds proved in the previous section
and Huisken's monotonicity formula to show that a given solution
converges under rescaling to the unique symmetric self-similarly
shrinking solution described in Theorem \ref{to self sim thm}.

Let us consider a network $M_t$ for $t\in [0,T)$ that contracts to a
point $x_0$ on the $x^1$-axis as $t\nearrow T$. Consider a sequence of
positive real numbers $\lambda_i \nearrow \infty$ satisfying
$\lambda_i^2T>1$. We rescale the network under a sequence of parabolic
dilations given by
\begin{equation}\label{dilations}
  M_{\tau}^i:=\lambda_i\left(M_{\lambda_i^{-2}\tau+T}-x_0\right).
\end{equation}
It is easy to see that for each $i$ the rescaled network $M_{\tau}^i$
satisfies the mean curvature flow equation in the time interval
$[-\lambda_i^2T,0)$. Notice that this parabolic rescaling also
preserves the symmetry, the gradient estimates and the $120^\circ$
condition. In particular, the rescaled networks are solutions
described by a graph and  the ratio estimate from Corollary \ref{din
  dex cor} $$ d_\text{in}\leq C d_\text{ex} $$ is preserved.

From the geometric estimates in the proof of Proposition
\ref{curvature bounds prop}, we deduce
\begin{corollary}\label{cinfinityestimates}
  The networks $M_\tau^i$ are uniformly bounded in $C^\infty$
  independently of $i$ on any compact time subinterval of $(-1,0)$.
\end{corollary}

\begin{proof}[Proof of Theorem \ref{self sim unique thm}]
  Note that since the networks $M_t$ (respectively $M_\tau^i$) satisfy
  the mean curvature flow equation they can be seen as Brakke flows
  with equality (see \cite{MantegazzaNetworks}).

  Huisken's monotonicity formula (\cite{Huisken90,MantegazzaNetworks})
  implies that the Gaussian density ratio given by
$$ \Theta (M_t,x_0,T) =
\int\limits_{M_t}\frac{1}{(4\pi(T-t))^\frac{1}{2}}
\exp\left(-\frac{|x-x_0|^2}{4(T-t)}\right)\,d\mu\qquad 
\hbox{for } t < T $$
is monotonically decreasing in time. Therefore, the limit
$$\Theta(x_0,T):=\lim_{t\nearrow T}\Theta (M_t,x_0,t)$$ 
exists and is finite. This limit is known as the Gaussian density and
it satisfies
$$\Theta(M_t,x_0,t)-\Theta(x_0,T)=
\int\limits_t^T\int\limits_{M_t}\rho_{x_0,T}(x,t)\,
\left|\kappa-{\frac{\langle
      x-x_0,\nu\rangle}{2(T-t)}}\right|^2d\mu_t(x)\,dt,$$ where
$\rho_{x_0,T}(x,t):=(4\pi(T-t))^{-1/2}\exp\left(-{\frac{|x-x_0|^2}{4
      (T-t)}}\right)$.

Changing variables according to the the parabolic rescaling described
by \eqref{dilations} we obtain
\begin{equation}
  \Theta(M_t,x_0,t)-\Theta(x_0,T)
=\int\limits_{\lambda_i^2(t-T)}^0\;\;\int\limits_{M_\tau^i}
\rho_{0,0}(y,\tau)\,\left|\kappa+{\frac{\langle y,\nu\rangle}{2\tau}}
\right|^2d\mu_\tau(y)\,d\tau.  
\label{densityidentity}
\end{equation}

Now for each $\lambda_i$ we look at the flow $M_\tau^i$ at time
$\tau=-1/2$. This defines a sequence $t_i \rightarrow T$ satisfying
$-\lambda_i^{-2}\frac12+T=t_i$. Using Corollary
\ref{cinfinityestimates} we can extract a subsequence, which we again
call $M_\tau^i$, that converges to a limit network $M^{\prime}_\tau$
and satisfies all the estimates above and, as before, flows by mean
curvature.

Moreover,  \eqref{densityidentity} implies
$$\int\limits_{-1}^0\;\int\limits_{M_\tau^i}
\rho_{0,0}(y,\tau)\,\left|\kappa+{\frac{\langle y,\nu\rangle}{2\tau}}
\right|^2d\mu_\tau(y)\,d\tau \rightarrow0.$$ 
Thus $M^{\prime}_{-\frac{1}{2}}$ satisfies
$$ \kappa = \langle x, \nu \rangle $$
and $(M^\prime_\tau)_{\tau < 0}$ is a self-similarly shrinking
solution which by the results of Section \ref{selfsimilar sec} and the
symmetry assumption is unique. \par Since the limiting flow does not
depend on the sequence $(\lambda_i)_i$ chosen, we obtain the stated
result by choosing an appropriate sequence $(\lambda_i)_i$.
\end{proof}

\section{Existence and Uniqueness of Self-Similar Lenses}
\label{selfsimilar sec}
In this section we prove Theorem \ref{to self sim thm}. We divide the
proof into two subsections: existence of self-similar shrinking lenses
and uniqueness in this class.

 \subsection{Existence}

We investigate the existence of a homothetically shrinking
solution $u(x,t)$ of \eqref{free bvp}.  That is, we study existence of
a profile $u_0$, a scaling
function $\lambda=\lambda(t)$  and an interval $I$ such that
\begin{equation*}
  \{ (x,u(x,t)): x \in \lambda(t)I \}= \lambda(t) \cdot  \{
  (y,u_0(y)): y \in I  \}. 
\end{equation*}
By letting $x=\lambda(t) y$ we obtain
\begin{equation}\label{homothetical solution}
u(x,t)=  \lambda(t) \cdot u_0\left( \frac{x}{\lambda(t)}\right).
\end{equation}
{}From  the evolution equation \eqref{free bvp}, we get
\begin{equation*}
\frac{u_0''}{1+(u_0')^2}=
\lambda(t)\dot{\lambda}(t)\left(-\frac{x}{\lambda(t)} u_0'+u_0\right)
=\lambda(t)\dot\lambda(t)(-yu_0'+u_0).
\end{equation*}
This implies $\dt\left(\tfrac12\lambda^2\right)=\text{const}$.  Since
$u(x,t)$ is a shrinking solution, it necessarily holds
$\dot\lambda(t)<0$, $\lambda(t)>0$. By assuming that $\lambda(-1/2)=1$
and fixing the blow up time to
be $T=0$ we obtain $\lambda(t)=\sqrt{-2t}$. For simplicity, in what
follows of this section we relabel $u:=u_0$ and $x:=y$. Then, we
obtain the following equation that determines the self-similarly
shrinking solution:
\begin{equation}\label{self sim}
u''=\left(1+u'^2\right)\left(xu'-u\right).
\end{equation}
(Note that this equation is equivalent to
$\kappa=\langle(x,u),\nu\rangle$.)  Thus, the existence statement of
Theorem \ref{to self sim thm} can be reduced to the following
proposition.

\begin{proposition}\label{self-similar solution}
  There exists a solutions $u$ to \eqref{self sim} that via
  \eqref{homothetical solution} induces a homothetically shrinking
  solution to \eqref{free bvp} for $-\infty<t<0$. Moreover, this
  solution is symmetric with respect to the $x^2$-axis.
\end{proposition}
\begin{proof}
Given $0<h<1$ let us consider the initial value problem
\begin{equation}\label{self sim symm}
\begin{cases}
u_{xx}= (1+u_x^2)(u_xx-u),&0\le x\le x_{\text{max}},\\
u(0)=h,&  \\ 
u_x(0)=0.&
\end{cases}
\end{equation}
We denote the solution to this equation by $u^h$.  In order to prove
the lemma we show that there is an $0<H<1$ and a point $x_0^H$ such
$u^H(x_0^H)=0$, $ u^H(x)>0$ for $x\in\left[0,x^h_0\right)$ and
$u^H_x(x_0^H)=-\sqrt{3}$. Then $u^H$ induces a homothetically
shrinking solution symmetric with respect to the $x^2$-axis.

In order to find $H$, consider
\begin{align*}\mathcal{S}:=\left\{\rule{0mm}{2.5ex}h\in(0,1):\right.& 
  \exists\, x^h_0\in\left(0,\sqrt2\right] \hbox{ such that }
  u^h(x)>0\hbox{ for }x\in\left[0,x^h_0\right),\\ & \left.  \
    u^h\left(x^h_0\right)=0\text{ and }
    \left|u^h_x\right|<\sqrt3\text{ on
    }\left[0,x^h_0\right]\right\}\end{align*} and define
$H:=\sup\mathcal S$. We will show that this $H$ satisfies the
properties required above.

In order to prove this we show that $\mathcal{S}$ is not empty and
open.  Moreover, we prove that $H<1$. This implies that for $h=H$
necessarily one of the conditions that define $\mathcal{S}$ has to be
violated. We will conclude the result from the existence of $x^H_0$
such that $ u^H(x)>0\hbox{ for }x\in\left[0,x^H_0\right),\
u^H\left(x^H_0\right)=0\text{ and }
\left|u^H_x\right|\left(x^H_0\right)\geq \sqrt3$. 

Lines through the origin are solutions to \eqref{self sim}. Hence
uniqueness of solutions to ordinary differential equations implies
that our solutions stay concave; i.e. $u^h_{xx}<0$. Since
$(xu^h_x-u^h)_x=xu^h_{xx}<0$ we get
\begin{equation*}
  u^h_{xx} =
  \left(1+\left(u^h_x\right)^2\right)\left(xu^h_x-u^h\right) \le
  xu^h_x-u^h \le -u^h(0)=-h,
\end{equation*}
and therefore $u^h(x) \leq h - \frac{1}{2}h x^2$. In particular, if
the gradient of $u^h$ remains bounded, there is an $0\leq x^h_0 \leq
\sqrt{2}$ such that $u^h(x)>0$ for $x\leq x^h_0$ and $u^h(x^h_0)=0$.

It is easy to check that $u^0\equiv0$.  Continuous dependence on the
initial data implies that for $|h|\ll1$ sufficiently small
$|u^h_x|<\frac{\sqrt3}2$ holds in $\left[0,\sqrt2\right]$. Thus we can
fix $\epsilon>0$ such that there exists $x_0^\epsilon>0$ satisfying
$u^\epsilon(x_0^\epsilon)=0$, and $u^\epsilon>0$,
$u^\epsilon_x>-\sqrt3$ in $[0,x_0^\epsilon]$. In particular, $
\mathcal{S}\ne \emptyset$.

On the other hand, it is also easy to check that $\graph u^1$ lies on
the unit circle. Hence, there is $x^1_{\sqrt3}>0$, $0<x^1_{\sqrt3}<1$,
such that $u^1_x(x^1_{\sqrt3})=-\sqrt3$ and $u^1>0$ in
$\left[0,x^1_{\sqrt3}\right]$.  Denote by $x^h_{\sqrt3}$ the value where
$$u^h_x\left(x^h_{\sqrt3}\right)=-\sqrt3.$$ 
As solutions are strictly concave, $x^h_{\sqrt3}$ depends smoothly on
$h$ near $h=1$.  Since $u^1$ does not extend up to $x=1$ and
$u^1(x)>0$ for $x\in [0,1)$, it holds that $H<1$.

Now, in order to show that $\mathcal{S}$ is open we note that $u^h$ is
strictly decreasing in $x$ for $x>0$ and $h>0$ and the value $x^h_0>0$
(where the solution $u^h$ becomes zero), depends smoothly on $h$ if
$u^h_x\left(x^h_0\right)>-\infty$. Similarly, since for $h\in
\mathcal{S}$ the solution $u^h_x$ does not blow up in $\left[0,
  x^h_0\right]$, $u^h_x$ depends smoothly on $h$ as well. This implies
$\mathcal{S}$ is open.

Consider $h_i\in \mathcal{S}$ such that $h_i\to H$. The previous
claims imply that $\lim x^{h_i}_0=x^H_0\le\sqrt2$ and
$u^H_x\left(x^H_0\right)\ge-\sqrt3$, which finishes the proof.
\end{proof}

\subsection{Uniqueness of  Self-Similar Lenses}

The main result of this section is:
\begin{lemma}\label{unique self-sim sol}
  The solution obtained in Proposition \ref{self-similar solution} is
  unique in the class of lens-shaped networks. In particular, the
  solution is symmetric with respect to the $x^2$-axis, that is
  $u_x(0)=0$.
\end{lemma}

We divide the estimates in this section in two types: geometric
estimates, that provide information for small heights $h$ (for the
definition of height see Section \ref{geomest}), and integral
estimates inspired by \cite{AndrewsClassification}.

\subsubsection{Geometric Estimates} \label{geomest} 

\begin{figure}[htb]
\epsfig{file=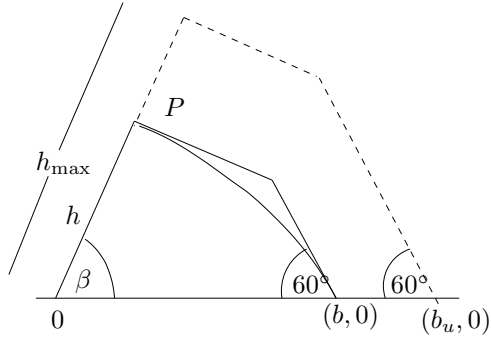, width=0.5\textwidth}
\caption{Circumscribed quadrilateral}
\label{quadrilateral fig}
\end{figure}

Let $P$ be the point where $\graph u$ attains the smallest distance to
the origin. In polar coordinates we describe this point by $r=h$ and
$\phi=\beta$.  Without loss of generality we assume that this point
lies on the right upper quadrant. By defining $b$ as the distance from
the origin to the point where $u_0$ intersects the $x^1$-axis we have
the following picture: see Figure \ref{quadrilateral fig}.

\psfrag{6}{$60^\circ$} 
  \psfrag{P}{$P$} 
  \psfrag{rmax}{$h_{\text{max}}$} 
  \psfrag{r}{$h$} 
  \psfrag{b}{$\beta$} 
  \psfrag{0}{$0$}
  \psfrag{a}{$(b,0)$}
 \psfrag{amax}{$(b_{u},0)$}
 \psfrag{rmax}{$h_{\text{max}}$}

In particular, lenses that are symmetric with respect to the $x^2$
axis satisfy $\beta=\frac{\pi}{2}$.

We use the following energy, which is a first integral for the
equation of a self-similarly shrinking solution:
\begin{equation}\label{energy 1}E^u:=\langle F,\nu\rangle e^{-\frac12|F|^2}
  \equiv\frac{u-xu_x}{\sqrt{1+u_x^2}}e^{-\frac12\left(x^2+u^2\right)}.
\end{equation} 
This energy was already applied by Uwe Abresch and Joel Langer
\cite{AbreschLanger} to classify all closed and immersed
self-similarly shrinking curves. Later on Ben
Andrews~\cite{AndrewsClassification} applied a modified energy to
classify such solutions when the normal speed of the evolution is a
positive power of the curvature.
\begin{lemma}\label{energy lemma}
  The energy $E^u$
  is constant in space along self-similar solutions $F$ of \eqref{mcf}
  with $T-t=\frac12$ or, equivalently, independent of $x$ for
  solutions $u$ to \eqref{self sim}.
\end{lemma}
\begin{proof}
  Direct computations using \eqref{self sim} give
  \begin{align*}
    E^u_x\sqrt{1+u_x^2}^3e^{\frac12\left(x^2+u^2\right)}=&\,
    -xu_{xx}\left(1+u_x^2\right)-u_xu_{xx}(u-xu_x)\\
    &\,-(u-xu_x)\left(1+u_x^2\right)(x+uu_x)\\
    =&\,u_{xx}\left[-x\left(1+u_x^2\right)-u_x(u-xu_x)+(x+uu_x)\right]=0.
  \end{align*}
\end{proof}

\noindent In what follows, when we use explicit numbers, the
corresponding estimates can  be verified easily by hand.

In what remains of this subsection we will limit ourselves to the case
$h\leq 0.5869\equiv h_{\text{max}}$.

Moreover, suppose that $(b,0)$ corresponds to the intersection point
of the self-similar curve and the positive $x^1$-axis. 
Then the following lemma holds:
\begin{lemma} If $h\leq h_{\text{max}}$ and a solution intersects the
  $x^1$-axis at an angle of $60^\circ$, then $ b\leq 0.7645\equiv
  b_{u}$ or $b\geq 1.2568\equiv b_{l}$
\end{lemma}
\begin{proof}
By Lemma \ref{energy lemma}, applied to $P$ and $(b,0)$, we have that
$$ he^{-\frac{h^2}{2}}=\sin \frac{\pi}{3}\cdot  be^{-\frac{b^2}{2}}.$$
Since $0\leq h\leq h_{\text{max}}$ and $h\mapsto h \exp(-h^2/2)$ is
increasing on $[0,1]$, we have
$$\sin \frac{\pi}{3}\cdot    be^{-\frac{b^2}{2}}\leq 0.5869\,
e^{-\frac{0.5869^2}{2}}\leq 0.49405.$$
Moreover,
\begin{align*} \sin \frac{\pi}{3}\cdot be^{-\frac{b^2}{2}}\geq&\
  \sin \frac{\pi}{3}\cdot 0.7645 e^{-\frac{0.7645^2}{2}}\ge0.49430
  \hbox{ for } 1 \geq b\geq 0.7645, \\
  \sin \frac{\pi}{3}\cdot be^{-\frac{b^2}{2}}\geq&\ \sin
  \frac{\pi}{3}\cdot 1.2568e^{-\frac{1.2568^2}{2}}\ge0.49408 \hbox{ for
  }1\leq b\leq 1.2568.
\end{align*}
\end{proof}

\noindent The following proposition states that if $h \leq
h_{\text{max}}$,  there are no solutions
to \eqref{self sim} such that $\graph u$ intersects the
$x^1$-axis on both sides at an angle of $60^\circ$ .

\begin{proposition}\label{small h prop}
  Assume that $0<h\le h_{\text{max}}$.  Then there is no solution $u$
  to \eqref{self sim} such that $u_x(b)=-\sqrt3$ and $u_x(a)=\sqrt3$,
  where $a$ is defined as $b$, but for the negative $x^1$-axis. 
\end{proposition}

The proof will be done by contradiction and it is divided into two
cases: $b\geq b_{l} $ and $b\leq b_{u}$. In order to obtain a
contradiction in the first case we will build a barrier function, for
which we will prove that it intersects the $x^1$-axis at $x^1<b_{l}$.
For the latter case we will consider the polygon depicted by doted
lines in Figure \ref{quadrilateral fig}. This polygon has fixed sides
$h_{\text{max}}$ and $b_{l}$, but changes with the angle $\beta$.
Moreover, if $b\leq b_{u}$, this polygon always contains the part of
the solution shown in the Figure \ref{quadrilateral fig}. We get a
contradiction by showing that the area of the polygon (denoted by
$A_\Box$) is always smaller that the one enclosed by the solution. In
order to show this we will need:

\begin{lemma}\label{der quad lem}
  Let $A_\Box$ be the area of the polygon with sides of length
  $h_{\text{max}}$ and $b_{u}$, respectively, as shown in Figure
  \ref{quadrilateral fig}.  Then we have
$$\frac{d A_\Box}{d\beta}<\frac{1}{2}.$$
\end{lemma}
\begin{proof}
  Observe that $\frac{dA_\Box}{d\beta}=\frac12l^2$, where $l$ denotes
  the length of the diagonal of the quadrangle that starts at the
  origin. This can be seen by bisecting the quadrangle along that
  diagonal and then regarding the variation of area as $\beta$
  increases as the extra triangle that appears next to this bisecting
  line. \par We conclude by noting that $l^2$ is estimated from above
  by $h_\text{max}^2+b_{u}^2<1$.
\end{proof}

Now we proceed with the proof of the Proposition.

\begin{proof}[Proof of Proposition \ref{small h prop}]

  We first want to remark that results in 
  \cite{AndrewsClassification,AbreschLanger}
  show that between points of minimum distance and maximum distance to
  the origin there is at least an angle of $\frac{\pi}{2}$. 
  As $h\le1$,
  the point $P$ has minimal distance to the origin. Hence the
  solutions stay graphical at least on that interval and it can be
  continued below the $x^1$-axis for an angle of
  $\frac{\pi}{2}-\beta$.  Now we proceed with the two cases described
  above.

\medskip

\noindent{\bf  Consider $b>b_{l}$ :}

\noindent
We first consider the ``symmetric case'', i.\,e.\ the case $u_x(0)=0$:
Assume that $x>0$.  We know that $u_x<0$ and $u_{xx}<0$. Moreover,
direct calculations and equation \eqref{self sim} imply that
$u_{xxx}<0$ and $u_{xxxx}<0$. By defining $f(x):=xu_x-u$ we obtain for
$k\in\N$
\[\frac{d^kf}{dx^k}(x)=(k-1)\frac{d^ku}{dx^k}(x)
+x\frac{d^{k+1}u}{dx^{k+1}}(x).\] Thus $f_x(x)=xu_{xx}<0$, $f_{xx}<0$,
and $f_{xxx}<0$. Hence
\begin{align*}
  (xu_x-u)(x)=&\,f(x)\le f(0)+f_x(0)\cdot x
  +\tfrac12f_{xx}(0)\cdot x^2<-h-\tfrac12hx^2,\\
  u_{xx}(x)=&\,\left(1+u_x^2\right)(xu_x-u)
  <-h\left(1+u_x^2\right)\left(1+\tfrac12x^2\right).
\end{align*}
As $u_x^2\ge0$, we get $u_x(x)<-hx-\tfrac16hx^3$ and
$u(x)<h-\frac12hx^2-\tfrac1{24}hx^4$.  Solving for this equation we
see that $u(x)<0$ for $x\geq \sqrt{-6+2\sqrt{15}}$. In particular,
$0<b<\sqrt{-6+2\sqrt{15}}\equiv B_1$.  Using $x=0$ and $x=b$, we get
according to Lemma \ref{energy lemma}
\[E^u=he^{-\frac12h^2}=\sin\frac{\pi}{3}\cdot b\cdot
e^{-\frac12b^2}.\] As $x\mapsto xe^{-\frac12x^2}$ is monotonically
increasing for $0\le x\le1$ and monotonically decreasing for $1\le x$,
we get for $h\le H_1\le1$
\begin{equation}\label{use Energy}
  H_1e^{-\frac12H_1^2}\ge E^u> \sin\frac{\pi}{3}\cdot B_1e^{-\frac12B_1^2}.
\end{equation}
A direct computation shows that this inequality is violated if
$0<h\le H_1:=0.5587$.

We assume now that $H_1\le h \leq h_{\text{max}}$. Then we have
$u_x^2\ge h^2\left(x+\tfrac16x^3\right)^2$. Thus,
\begin{align*}
  u_{xx}(x)\le&\,-h\left(1+h^2\left(x+\tfrac16x^3\right)^2\right)
  \left(1+\tfrac12x^2\right)\\
  =&\,-h\left(1+\left(\tfrac12+h^2\right)x^2+\tfrac56h^2x^4
    +\tfrac{7}{36}h^2x^6+\tfrac{1}{72}h^2x^8\right),\\
  u_x(x)\le&\,-h\left(x+\tfrac13\left(\tfrac12+h^2\right)x^3
    +\tfrac1{6}h^2x^5+\tfrac1{36}h^2x^7+\tfrac1{648}h^2x^9\right),\\
  u(x)\le&\,h-h\left(\tfrac12x^2
    +\tfrac1{12}\left(\tfrac12+h^2\right)x^4
    +\tfrac1{36}h^2x^6+\tfrac{1}{288}h^2x^8
    +\tfrac1{6480}h^2x^{10}\right).
\end{align*}
We obtain $b<1.2568$. Arguing as above, this yields
$h>h_{\text{max}}$, which contradicts the assumptions of the lemma.

Now we consider the case $u_x(0)\neq0$: We may assume that $u_x(0)<0$.
We start by rotating the coordinates such that the new $x^2$-axis
agrees with the old line of slope $\tan \beta$. Using the first
barrier from the argument in the symmetric case we see that we get a
contradiction if the $x^1$-coordinate of the intersection of $\graph
u$ and the new $x^1$-axis is bigger than $B_1$.  Results in
\cite{AndrewsClassification} imply that the distance to the origin is
increasing up to the point where the maximum is attained. Thus we also
get a contradiction if the distance to the origin at the intersection
with the old $x^1$-axis is bigger than $B_1$. Since we have assumed
that our solution intersects this old $x^1$-axis at an angle of
$60^\circ$, we can exclude initial heights $h$ with $0<h\le H_1$ as
before. Repeating this argument with the second barrier from the
symmetric case yields the stated claim.

\medskip

\noindent{\bf Consider $b< b_{u}$:} 

\noindent As mentioned above, in this situation the strategy of the
proof is to compute and compare the difference between the areas
enclosed by the quadrilateral shown in Figure \ref{quadrilateral fig}
and by the solution. We will observe that for $0<h\leq
h_{\text{max}}$ this difference is always negative, contradicting the
definition of the polygon.

From \cite{AndrewsClassification}, we obtain that the angle between a
minimum and a maximum of $|F|=\sqrt{x^2+u^2}$ is at most
$\frac\pi{\sqrt 2}$.  Moreover, this angle is at least $\frac\pi2$.
Thus the boundary conditions imply that there is precisely one minimum
of $|F|$ above the $x^1$-axis.  Hence $\beta>\pi-\frac\pi{\sqrt2}$,
for otherwise, $|F|$ would have a local maximum between $P$ and the
negative $x^1$-axis.

As before, the area $A$ of the domain bounded by the self-similarly
contracting solution and the half-lines corresponding to $\theta=0$
and $\theta=\beta$ can be computed using the divergence theorem and
the fact that $\kappa = \langle(x,u), \nu \rangle$. We
obtain
\[A=\frac12\int\limits_\gamma\kappa\,d\mu =\tfrac12\left(2\pi
  -\left(\pi-\frac\pi3\right) -\left(\pi-\beta\right)
  -\left(\pi-\frac\pi2\right)\right) =\frac\beta2-\frac\pi{12},\]
where $\gamma$ denotes the curved part of the boundary of the domain
and $d\mu$ the corresponding volume element.

Recall from Figure \ref{quadrilateral fig} that for every fixed
$\beta$, $b\leq b_{u}$ and $ h\leq h_\text{max}$ there is a natural
common polygon (in doted lines) bounding from above the area of the
self-similar curve. A simple computation shows that the area of the
quadrilateral $A_\Box(\beta)$ is given by
\begin{gather} \begin{split}\label{area quad}
A_\Box(\beta)=&\ \frac{b_{u}^2\cos\beta \sin \beta}{2}+
(h_\text{max}-b_{u}\cos\beta) b_{u}\sin\beta
-\frac{(h_\text{max}-b_{u}\cos\beta)^2\tan \alpha}{2}\\ =&\
\frac{b_{u}^2}{4} \sin 2\beta+(h_\text{max}-b_{u}\cos\beta)
b_{u}\sin\beta
-\frac{(h_\text{max}-b_{u}\cos\beta)^2\tan \alpha}{2},
\end{split} \end{gather}
where $\alpha=\pi-\frac{\pi}{3}-\beta=\frac{2\pi}{3}-\beta$.
Furthermore, Proposition \ref{der quad lem}  implies that
$$\frac{d A_\Box}{d\beta}-\frac{d A}{d\beta}<0.$$
Hence an explicit calculation shows that 
$$(A_\Box-A)(\beta)\leq 
(A_\Box-A)\left(\pi-\tfrac\pi{\sqrt2}\right)
<0.$$
\end{proof}

\subsubsection{Integral Estimates}
In order to finish the proof of uniqueness we use the support
function. Recall that this is defined to be the map $S:S^1\times [0,T)
\rightarrow \mathbb{R}$ given by
\[S(x,t)=\left\langle x,F(\nu^{-1}(x),t) \right\rangle,\] where $\nu$
is the function $G\rightarrow S^1$ giving the unit normal vector to
the curve. Note that $F(\nu^{-1}(x),t)$ is the point of the curve
where $x$ is the unit normal vector. Furthermore $S^1$ is parameterized
over its arc-length $\theta$, such that we get $S$ as a function of
$(\theta,t)\in [0,2\pi]\times [0,T)$.  In general, the support
function and the curvature are related by the formula
\[\kappa^{-1}=S_{\theta \theta}+ S.\] Hence, for self--similar
shrinking solutions one has
\[S(\theta)=\kappa(\theta).\] and the equation for the support
function reads: \begin{equation} \label{eq supp} S_{\theta \theta}+
  S=\frac{1}{S}. \end{equation} In the following we will assume that
$S>0$. An easy calculation shows that the
energy $E^S$ defined as
\begin{equation} \label{def energy}
E^S(S,S_{\theta}):= (S_{\theta})^2+ S^2 -2\log S
\end{equation}
is a first integral of the Equation \eqref{eq supp}. Therefore, each
solution $S(\theta)$ of equation \eqref{eq supp} lies on a level set
of $E^S$. Note that this energy is equivalent to the energy considered
in Lemma \ref{energy lemma} since $-2\log E^u=E^S$. 

Let us have a closer look at the energy levels in the
$(S,S_{\theta})$--plane: $E^S$ has a unique critical point at $(1,0)$.
The fact that \[D^2E^S(S,R)=\left( \begin{array}{cc} 2+ \frac{2}{S^2} &
    0 \\ 0 & 2
\end{array}  \right) \] is positive definite implies that the other
level sets are convex closed curves around $(1,0)$. See also the
energy diagram picture. \\ 
\psfrag{SSS}{$S$}
\psfrag{Stheta}{$S_{\theta}$}
\psfrag{S1}{$S_1$}
\psfrag{O}{$(0,1)$}
\psfrag{S2}{$S_2$}
\psfrag{E}{$E^S$}
\psfrag{S-}{$S_-$}
\psfrag{S+}{$S_+$}
\begin{figure}
\includegraphics[width=0.5\textwidth]{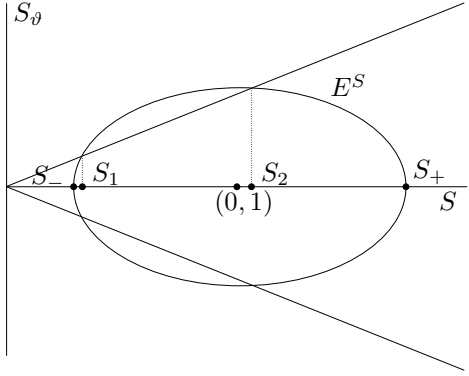}
\caption{Energy diagram.} 
\end{figure}

Furthermore, it holds \begin{itemize} 
\item $E^S(S,R)>0$ for all $S>0$,  
\item $E^S(1,0)=1$,
\item $E^S(S,0) \longrightarrow \infty$ for $S \longrightarrow 0$ or $S
  \longrightarrow \infty$. 
\end{itemize}
In order to impose the angle condition at the triple points we
compute the points of each energy level where the embedding vector $X$
and the tangent vector $X_{\theta}$ form an angle of $\pm \frac{2
  \pi}{3}$ or $\pm\frac\pi3$. (Later we will consider integrals that
start at a point corresponding to a minimum of $|F|$. Thus we exclude
the undesired angles.) This condition is equivalent to
\begin{equation} \label{2lines} \frac12=\abs{\cos
    \frac{\pi}{3}}=\abs{\cos \frac{2\pi}{3}}=\frac{\abs{\langle
      X,X_{\theta}\rangle}}{\abs{X}\abs{X_{\theta}}}.
\end{equation}
Now let us express $X$ and  $X_{\theta}$ in terms of  $S$ and
$S_{\theta}$: 
\begin{align*}
X=&S   \nu + S_{\theta}\nu_{\theta},\\
X_{\theta}=&S   \nu_{\theta} + S_{\theta \theta}\nu_{\theta}.
\end{align*}
Plugging these formulas in  \eqref{2lines} leads to \[
\frac{1}{2}=\frac{|S_{\theta}|}{\sqrt{S^2+ S_{\theta}^2}}, \] 
which is equivalent to 
\begin{equation} \label{22lines}
 S=\pm \sqrt{3} S_{\theta}.
\end{equation}

Therefore, in terms of the energy diagram above a self-similar lens
can be described as follows: It belongs to some energy level $E^S$ and
the points corresponding to $(a(t),0)$ and $(b(t),0)$ lie on the lines
$S=\sqrt3 S_\theta$ and $S=-\sqrt3 S_\theta$, respectively. Moreover,
these endpoints correspond to some values $\theta_a$ and $\theta_b$ in
the parameterization for the support function that satisfy
$\theta_a-\theta_b=\frac{2\pi}3$.

In the coming proof we make use of an analog of Ben Andrews' formula
\cite{AndrewsClassification} that allows us to compute the total
curvature between intersection points of a given energy level $E^S$
with the two half-lines $S= \pm \sqrt3 S_\theta$. In particular we
prove that the fact that the total curvature along the self-similar
shrinking curve should be $\frac{2\pi}3$ uniquely selects $E^S$ and
the lens. \\ For a fixed energy level $E^S$, let us define $S_-$ and
$S_+$ be the left and right intersections of the energy level $E^S$
with the $S$--axis.  Let $S_1$ and $S_2$ be the abscissas of the two
intersections of that energy level of $E^S$ with the lines
$S=\pm\sqrt3 S_\theta$, if they exist (see also the energy diagram).
Then, the integral of the curvature of the self-similar solution from
$S_-$ to $S_{1,2}$ is given by
$$\Psi\left(S_-,S_{1,2},E^S\right)
=\int_{S_-}^{S_{1,2}}\frac{1}{\sqrt{E^S-S^2+2\log 
    S}}\, dS.$$ For $S_1\le S_2$ the quantity \[
\tilde\eta=\frac{S_2}{S_-}\] parameterizes the energy levels (or at
least those for which we have intersections with the half-lines, which
are of interest). We will see that this also applies for $\bar
\eta=\frac{S_1}{S_-}$. Set $x=\frac{S}{S_-}$ and
\begin{equation}\label{psi(eta)}
  \Psi(\eta)=\int\limits_{1}^{\eta}\frac{1}{\sqrt{ 1-x^2+
      \frac{4\eta^2-3}{3\log \eta} \log x}}\, dx. 
\end{equation} 
Thus, the total curvature of a symmetric lens is given by
$2\Psi(\eta)$. Similarly, for asymmetric lenses the total curvature
can be computed by $\Sigma(\eta)=\Psi(\eta)+\Psi(\bar{\eta})$.
(According to \cite{AndrewsClassification} none of our self-similarly
shrinking lenses contains a point corresponding to $S_+$. This ensures
in particular that the intersection angles are as desired.) 
\par From the energy conservation we notice that (for a fixed lens)
$\tilde\eta$ and $\bar{\eta}$, defined as above, satisfy the following
relation
\begin{equation}\label{def bar eta}\frac{4 \ol\eta^2 -3}{3 \log
    \ol \eta}=\frac{4 \tilde\eta^2 -3}{3 \log\tilde\eta}. \end{equation}
Moreover,
Lemma \ref{eta estimates} \eqref{A increasing} implies the following:
For each $\tilde\eta$ there is at most one $\bar{\eta}\ne \tilde\eta$ 
such that \eqref{def bar eta} is satisfied and, thus,  we can consider 
$\bar{\eta}$ as a function of $\tilde\eta$.
Moreover, there is a unique $\eta_0$ such that 
$\bar{\eta}(\eta_0)=\eta_0$ and  the function $\bar{\eta}(\tilde\eta)$ 
is decreasing for $\tilde \eta\geq \eta_0$. Hence, $\bar{\eta}\leq \eta_0$
can also be considered as a function of
$\tilde\eta\geq \eta_0$.

In this section we will use the following notation conventions: When
we use the variable $\eta$, we refer to $\eta\in[1,\infty)$. If we
restrict to values of $\eta$ to $\eta\geq\eta_0$, we will denote the
variable by $\tilde \eta$.  We will always assume that $1\le
\bar{\eta}\leq \eta_0\leq \tilde\eta$ and we consider $\bar{\eta}$ as
a function of $\tilde \eta$.

We want to remark, that the above considerations imply that $\Sigma$
is a function of $\eta$ and that it is enough to study
$\Sigma(\tilde \eta)$ where $\tilde \eta\geq \eta_0$.

Since the total curvature of a self-similarly shrinking solution is
$\frac{2\pi}{3}$, our uniqueness result follows from

\begin{proposition}\label{pi drei una}\neueZeile
  \begin{enumerate}[(i)]
  \item \label{psi prop} The function $\Psi(\eta)$ attains the value
    $\frac\pi3$ at a unique point $\eta \in[1,\infty)$.
  \item \label{sigma prop}$\Sigma(\tilde\eta)< \frac{2\pi}{3}$ for every
    $\tilde\eta\ge \eta_0$.
  \end{enumerate}
\end{proposition}

Moreover, we can draw the following connection with the previous
subsection:

\begin{lemma}\label{1.9 fine}
  For $\eta\ge1.9$ we have $\Psi(\eta)\neq\frac\pi3$ and
  $\Sigma(\tilde\eta)\neq\frac{2\pi}3$. 
\end{lemma}
\begin{proof}
  Assume that $\Psi(\eta)=\frac\pi3$ or
  $\Sigma(\tilde\eta)=\frac{2\pi}3$.  Then this $\eta$ (or
  correspondingly $\tilde \eta$) determines a self-similarly shrinking
  solution with correct intersection angles. Recall that $\tilde \eta$
  denotes the same quantity as $\eta$ does, but restricted to
  $[\eta_0,\infty)$. Hence in the computations of this proof we will
  always refer to $\eta$, but it can be directly replaced by $\tilde
  \eta$ to obtain the desired conclusion for $\Sigma(\tilde \eta)$.
  
  By definition of
  $h$ we have $S_-=h$ and $\eta=\frac{S_1}{h}$. At the point where our
  solution makes an angle of $\pi/3$ with the ray from the origin, we
  get $S_\theta = \frac{1}{\sqrt{3}}S=\frac1{\sqrt3}S_1$. As the
  energy is conserved we have
  \[h^2-2\log h=E^S=S_\theta^2+S^2-2\log S=\tfrac43S_1^2-2\log S_1.\]
  Rearranging terms gives
  \[2\log\eta=\left(\tfrac43\eta^2-1\right)h^2.\] From Proposition
  \ref{small h prop}, it suffices to consider $h\ge h_{\text{max}}$.
  Hence
  \begin{equation}\label{eta h desigualdad}
    2\log\eta\ge\left(\tfrac43\eta^2-1\right)h_{\text{max}}^2.
  \end{equation}
  It is easy to check that for $\eta\ge1.9$ inequality \eqref{eta h
    desigualdad} is violated.
\end{proof}

The previous Lemma implies that in order to prove Proposition \ref{pi
  drei una} it suffices to focus on $1<\eta\leq 1.9$. In this range we
perform a careful analysis of the functions $\Psi$, $\Sigma$ and their
derivatives. We first compute
\begin{align}\label{psi'}
  \frac{d\Psi}{d\eta}(\eta)=&\,\frac{\sqrt{3}}{\eta}-A(\eta)
  \int\limits_1^\eta \frac{\log x}{\left(1-x^2+\frac{4 \eta^2 -3}{3 \log
        \eta} \log x\right)^{\frac 3 2}} dx 
\end{align}
with
$$  A(\eta):=\frac{8 \eta^2 \log \eta - 4 \eta^2 + 3}{6 \eta \log^2
  \eta}.$$ Thus, $\Psi$ is increasing when $A(\eta)\leq 0$. In order
to relate this computation to the analogous one for the function
$\Sigma$, from \eqref{def bar eta} we notice
\begin{equation}\label{der bar eta} A(\tilde\eta)=A(\ol \eta)\frac{d\ol
    \eta}{d\tilde\eta}. \end{equation} 
Prior to the proof of Proposition \ref{pi drei una} we need a few
elementary lemmata.

\begin{lemma}\label{eta estimates}
The following estimates hold
\begin{enumerate}[(i)] 
\item \label{A increasing} $\frac{dA}{d\eta}(\eta)\equiv A'(\eta)>0$
  for all $\eta>1$. In particular, there is a unique
  $\eta_0\in[1,\infty)$ such that $A<0$ on $[1,\eta_0)$ and $A>0$ on
  $(\eta_0,\infty)$. Notice that this agrees with the previous
  definition of $\eta_0$.
\item \label{33a lem} $ \frac{\sqrt{3}\eta}{3A(\eta)\log \eta}$ is
  decreasing for $\tilde\eta=\eta>1$ such that $A(\eta)>0$ and
  $\eta\leq e$ (in particular for $\eta\leq1.9$).
\item \label{monot in eta} The function $\eta \mapsto
  \left[1+\fracd{4\eta_0^2-3}{4\eta^2-3}\right]\log\eta$ is
  increasing for $\eta \geq 1$.
\end{enumerate}
\end{lemma}
\begin{proof}\neueZeile
\begin{enumerate}[(i)]
\item Define
\begin{equation*}
   B(\eta):=\left( 6 \eta^2 \log^3 \eta \right)
            \frac{dA}{d\eta}(\eta)=8 \eta^2 \log^2 \eta -12 \eta^2
            \log \eta - 3 \log \eta + 8 \eta^2 - 6. 
\end{equation*}
We show that $B$ stays positive for all $\eta \geq 1$. This follows from
$B(1)=2>0$ and $\frac{dB}{d\eta}\geq 0$: 
\begin{gather*}
\begin{split}
  \frac{dB}{d\eta}(\eta)
  =&\, 16 \eta \log^2 \eta - 8 \eta \log \eta + 4 \eta - \frac 3 \eta\\
  =&\, \frac 1 \eta \left( (4\eta\log\eta-\eta)^2
    +3(\eta^2-1)\right)\geq 0.
\end{split}
\end{gather*}

\item  By direct calculation, we
  obtain
\begin{align*}
  \frac{d}{d\eta}\frac{\sqrt{3}\eta}{3A(\eta)\log \eta}=&\,
  -\frac{\sqrt{3}A'(\eta)\eta}{3A^2(\eta)\log \eta}
  +\frac{\sqrt{3}}{3A(\eta)\log^2\eta} \left(\log \eta-1\right).
\end{align*}
 From \eqref{A increasing} the first term is negative.
Since $A(\eta)>0$, the second term is negative for $\log \eta-1\leq 0$
or equivalently $\eta\leq e$.

\item 
Set $c=4\eta_0^2-3$.
Taking a derivative we have
\begin{align*}\frac{d}{d\eta} \left(\left[1+
      \frac{c}{4\eta^2-3}\right] \log\eta\right)
  =&\ \frac{\eta}{(4\eta^2-3)^2}\left[\left(4\eta
      -\frac3\eta\right)^2+4c
    -\frac{3c}{\eta^2}-8c\log \eta\right]\\ 
  \ge&\ \frac{\eta}{(4\eta^2-3)^2}\left[\left(4\eta
      -\frac3\eta\right)^2+16
    -\frac{15}{\eta^2}-20\log \eta^2\right].\\ 
\end{align*}
Here we have used that $A(1.33652)>0$,\ $A(1.3365)<0$, hence
$1.33652>\eta_0> 1.3365$ and $5>c>4$. Taking derivatives, a direct
computation implies that the last square bracket above is increasing.
Since at $\eta=1$ this bracket is positive, we conclude the result.
\end{enumerate}
\end{proof}

\begin{lemma}\label{x estimates}\neueZeile
\begin{enumerate}[(i)]
\item \label{technical lemma} $$\frac{\log
    x}{1-x^2+\frac{4\eta^2-3}{3\log\eta}\log x} \leq
  \frac{3\log\eta}{\eta^2}$$ for all $1\leq x \leq \eta$ with strict
  inequality for $x<\eta$.
\item \label{inc quan} The function
\[ \frac{x-1}{1-x^2+\frac{4\eta^2-3}{3\log\eta}\log x}\equiv
  \frac{x-1}{d(x)}\]
is an increasing function of $x$, for $\eta\ge x\geq 1$.
\end{enumerate}
\end{lemma}
\begin{proof}\neueZeile
\begin{enumerate}[(i)]
\item The statement of the lemma is equivalent to $\frac{x^2-1}{\log
    x} \leq \frac{\eta^2-1}{\log \eta}$ for all $1\leq x \leq \eta$.
  Applying l'H\^opital's rule we see that $\frac{x^2-1}{\log
    x}\equiv\frac{f(x)}{g(x)}$ goes to $2$ for $x \to 1$. The
  functions $x^2-1$ and $\log x$ are both positive for $x>1$. Consider
  $f'g-g'f$. It vanishes at $x=1$ and is positive for $x>0$ as $x^2-1$
  is convex and $\log x$ is concave. Therefore the quotient is
  strictly increasing.
\item Notice that
\begin{align*}
 \frac{x-1}{d(x)}=&\
  \frac{1}{-(1+x)+\frac{4\eta^2-3}{3\log\eta}\frac{\log x}{x-1}}
\intertext{and}
\frac{d}{dx}\left(-(1+x)+\frac{4\eta^2-3}{3\log\eta}
    \frac{\log x}{x-1}\right)=&\
  -1+\frac{4\eta^2-3}{3\log\eta} \left(\frac{1}{x(x-1)}
    -\frac{\log x}{(x-1)^2}\right)\\  
  =&\ -1+\frac{4\eta^2-3}{3\log\eta} \frac{1}{(x-1)^2}
  \left(1-\frac{1}{x}-\log x\right)\\ 
  \leq&\
  \frac{4\eta^2-3}{3\log\eta}\frac{1}{(x-1)^2}
  \left(1-\frac{1}{x}-\log x\right). 
\intertext{Since}
\frac{d}{dx}\left(1-\frac{1}{x}-\log x\right)=&\
\frac{1}{x^2}-\frac{1}{x}\leq 0\hbox{ for }x>1,
\intertext{we conclude}
\frac{d}{dx}\left(-(1+x)+\frac{4\eta^2-3}{3\log\eta}
\frac{\log x}{x-1}\right)\leq&\ \frac{4\eta^2-3}{3\log\eta}
\frac{1}{(x-1)^2}\left(1-\frac{1}{1}-\log(1)\right)=0.
\end{align*}
\end{enumerate}
\end{proof}

\begin{proof}[Proof of Proposition \ref{pi drei una}:]
We prove \eqref{psi prop} and \eqref{sigma prop} of the  statement of 
Proposition  \ref{pi drei una} separately, intermitted by a lemma.
   
\noindent {\bf Proof of \eqref{psi prop}:}
We divide the proof into the following steps.
  \begin{enumerate}[(a)]
  \item\label{pi drei una i} We start by observing
    $\frac{d\Psi}{d\eta}>0$ for $1<\eta\le\eta_0$.
  \item\label{pi drei una ii} We show $\Psi>\frac\pi3$ at critical
    points of $\Psi$ in the interval $[1,1.9)$.
  \item\label{pi drei una iii} We prove $\Psi(1.9)>\frac\pi3$.
  \item\label{pi drei una iv} We conclude the result from the previous
    steps.
  \end{enumerate}
  
  Claim \eqref{pi drei una i} follows directly from \eqref{psi'} and
  Lemma \ref{eta estimates} \eqref{A increasing}.

For claim \eqref{pi drei una ii}, it suffices to consider
$\eta>\eta_0$ as $\Psi$ is strictly increasing otherwise.  Using Lemma
\ref{x estimates} \eqref{technical lemma}, we estimate
\begin{align}
  \begin{split}\label{dpsi}
    \frac{d\Psi}{d\eta}(\eta)&\,>
    \frac{\sqrt{3}}{\eta}-A(\eta)\frac{3\log \eta}{\eta^2} \int\limits_1^\eta
    \frac{1}{\left(1-x^2+\frac{4 \eta^2 -3}{3 \log
          \eta} \log x\right)^{\frac 1 2}} dx\\
    &\,=\frac{\sqrt{3}}{\eta}-A(\eta)\frac{3\log
      \eta}{\eta^2}\Psi(\eta).
  \end{split}
\end{align}

Suppose that $\Psi$ has a critical point for some
$\eta\in[1,1.9)\cap (\eta_0,1.9)$. Then $
\frac{d\Psi}{d\eta}(\eta)=0$ and thus
  $$\Psi(\eta)>\frac{\sqrt{3}\eta}{3A(\eta)\log \eta}.$$
  Lemma \ref{eta estimates} \eqref{33a lem} implies that
$$\Psi(\eta)>\frac{\sqrt{3}\cdot1.9}{3A(1.9)\log 1.9}\geq
\frac{\pi}{3},$$ proving \eqref{pi drei una ii}.

In order to prove claim \eqref{pi drei una iii}, we use claim
\eqref{pi drei una ii} and \eqref{dpsi}. According to Lemma \ref{1.9
  fine} and Proposition \ref{self-similar solution}, there exists
$\eta<1.9$, such that $\Psi(\eta)=\frac\pi3$.  Suppose that
$\Psi(1.9)\le\frac\pi3$. Then from \eqref{pi drei una ii},
$\frac{d\Psi}{d\eta}(1.9)\leq0$. However, from Equation \eqref{dpsi}, we
get $\frac{d\Psi}{d\eta}(1.9)>0$ which yields a contradiction.

Combining the previous steps, Proposition \ref{self-similar solution}
and Lemma \ref{1.9 fine}, \eqref{pi drei una iv} follows. \\[2ex]
We need the following useful lemma for the proof of \eqref{sigma prop}.
\begin{lemma} \label{bd for psieta0}
$$0.72<\Psi(\eta_0)\leq 0.785.$$
\end{lemma}

\begin{proof}
  It is easy to check that $A(1.33652)>0$ , therefore $\eta_0\leq
  1.33652$. Moreover, if $\Psi(1.33652)\leq \frac{\pi}{3}$, the proof
  of Proposition \ref{pi drei una} \eqref{psi prop} implies that
  $\Psi(\eta)$ is increasing for $\eta \leq 1.33652$ and
  $\Psi(\eta_0)\leq \Psi(1.33652)$. Hence, it is enough to show that
  $$\Psi(1.33652)\leq 0.785<\frac{\pi}{3}.$$ Set $\delta
  x:=\frac{1.33652-1}3$. Using Lemma \ref{x estimates} \eqref{inc
    quan} we have
  \begin{align*} \Psi(1.33652) \le&\,\int\limits_1^{1+3\delta
      x}\sqrt{\frac{x-1}{d(x)}}
    \sqrt{\frac1{x-1}}\,dx\umbruch\\
    \le&\,\sum\limits_{i=0}^{2}\int\limits_{1+i\delta
      x}^{1+(i+1)\delta x} \sqrt{\frac{(i+1)\delta x}{d(1+(i+1)\delta
        x)}} \sqrt{\frac1{x-1}}\,dx\umbruch\\
    =&\,\sum\limits_{i=0}^{2} \sqrt{\frac{(i+1)\delta
        x}{d(1+(i+1)\delta x)}}\ 2
    \left(\sqrt{(i+1)\delta x} -\sqrt{i\delta x}\right)\umbruch\\
    \le&\,0.785.  
  \end{align*}

  The lower bound can be computed analogously by noticing that
  $A(1.3365)<0$ and, thus, $\eta_0> 1.3365$. We remark that as $x\searrow
  1$ we consider the corresponding limit.
\end{proof}

\noindent{\bf Proof of \eqref{sigma prop}:}

Let $$f(\tilde\eta)=\Sigma(\tilde\eta)-\sqrt{3}\log(\tilde\eta \bar{\eta}).$$ 
Using \eqref{der bar eta}, we see that
\begin{align*}\frac{d f}{d \tilde\eta}=&\  
  \frac{d\Sigma}{d\tilde\eta}-\frac{\sqrt{3}}{\tilde\eta}-\frac{\sqrt{3}}
  {\bar{\eta}}\frac{d\bar{\eta}}{d\tilde\eta}\\ 
  =&\ -A(\tilde\eta)\left(
    \int\limits_{1}^{\tilde\eta}\frac{\log x}{\left( 1-x^2+
        \frac{4\tilde\eta^2-3}{3\log \tilde\eta} \log
        x\right)^{\frac{3}{2}}}dx \right. \\
  & \qquad\qquad\quad\left.+
    \int\limits_{1}^{\bar{\eta}}\frac{\log x}{\left( 1-x^2+
        \frac{4\tilde\eta^2-3}{3\log \tilde\eta} \log
        x\right)^{\frac{3}{2}}}dx \right) \leq 0
\end{align*}
as $\tilde\eta\geq \eta_0\ge\bar\eta$. Therefore, 
$$\Sigma(\tilde\eta)\leq\Sigma(\eta_0)-\sqrt{3}\log(\eta_0^2)+
\sqrt{3}\log(\tilde\eta\bar{\eta}).$$ Using Lemma \ref{bd for psieta0}
and that $\eta_0> 1.3365$ we get $$\Sigma(\tilde\eta)\leq 2\times
0.785-\sqrt{3}\log(1.3365^2)+ \sqrt{3}\log(\tilde\eta \bar{\eta}).$$
From \eqref{def bar eta} observe that
$$\log\ol\eta=\frac{4\ol\eta^2-3}{4\tilde\eta^2-3}\log\tilde\eta
  \le\frac{4\eta_0^2-3}{4\tilde\eta^2-3}\log\tilde\eta\ .$$
Hence, Lemma \ref{eta estimates} \eqref{monot in eta} implies for
$\eta_0\leq \tilde\eta \leq  1.9$ that 
\begin{align*} 
  \Sigma(\tilde\eta)\le&\,
\sqrt3\left[1+\frac{4\eta_0^2-3}{4\tilde\eta^2-3}\right]\log\tilde\eta
+  0.5653
\umbruch\\
      \le&\,\sqrt3\left[1+\frac{4(1.33652)^2-3}{4(1.9)^2-3}\right]
      \log(1.9)+ 0.5653
\umbruch\\
      <&\,\frac{2\pi}3.
    \end{align*}
\end{proof}

\begin{appendix}
\section{Existence and Uniqueness of Homothetically Shrinking
    Fish-Shaped Networks}\label{fish app}
\noindent
We aim to prove the existence and uniqueness of a fish-shaped
self-similarly shrinking network. By uniqueness we mean that
up to rotations and reflections,
this is the only self-similarly shrinking solution with two half-lines
going to infinity, where the half-lines are not parallel and the
network is topologically a lens.\\[2ex]

We generalize the class of networks we are working in. The name of the
following class of networks is inspired by Theorem \ref{fish thm} and
Figure \ref{fish pic}.
\begin{definition}[Generalized lens-shaped networks]
  A generalized lens-shaped network is an embedding $F:G\to
  M\subset\R^2$, where $G$ is an abstract graph as in the definition
  of a lens-shaped network. We impose the following conditions on $F$:
  \begin{enumerate}[(i)]
  \item $F$ is a homeomorphism on its image and, restricted to each
    edge, a diffeomorphism.
  \item At the vertices the images of the edges meet at $120^\circ$.
  \item As approaching infinity the non-compact edges are close to
    straight half-lines.
  \end{enumerate}
\end{definition}

\noindent Assume that a self-similarly shrinking network contains a
curve which is close to a half-line as approaching infinity. Then it
is easy to see that this part of the network has to be a straight
line.

\begin{theorem}\label{fish thm}
  There exists a family $(M_t)_{t\in(-\infty,0)}$ of generalized
  lens-shaped networks that has the following properties:
  \begin{enumerate}[(i)]
  \item The networks $(M_t)_t$ shrink homothetically under \eqref{mcf}
    away from the triple points.
  \item The bounded components of $\left(\R^2\setminus M_t\right)_t$
    shrink to the origin in Hausdorff distance as $t\nearrow0$.
  \item The family $(M_t)_t$ is different from the one in Theorem
    \ref{to self sim thm}. 
  \end{enumerate}
  Moreover, up to rotations, $(M_t)_t$ is unique.
  We call such a network $M_t$ fish-shaped. 
\end{theorem}

\begin{proof}
  We distinguish whether the network is reflection symmetric with
  respect to an axis or not. Notice that results of J\"org
  H\"attenschweiler \cite[Lemma 3.17]{HaettenschweilerDiplom} imply
  that the network has such a symmetry. However, for the reader's
  convenience, we give an independent non-existence proof in the
  ``asymmetric case''.
  \par \noindent {\bf Existence and Uniqueness in the Symmetric Case:}
  Note that the total curvature $K$ of the loop is given by
$$K=\int_{\text{loop}}\kappa\, ds = 2\pi-\frac{2}{3}\pi =
\frac{4\pi}{3}\ .$$ Note furthermore that the energy as defined in
\eqref{energy 1} is constant along the loop. At triple points this
follows from the fact that at each triple point there is a straight
half-line about which the network is locally symmetric.  Thus both
arcs which form the loop are part of the same self-similarly shrinking
solution (without triple junctions). From the angle condition at the
triple points there has to be at least one point of least distance
$r_\text{min}$ to the origin on each arc. As curved parts locally
solve \eqref{self sim}, \cite{AndrewsClassification} implies that
local and global extrema of the distance to the origin coincide. Let
us first assume that both arcs have only one point of least distance
to the origin (and thus no point of maximal distance $r_\text{max}$ to
the origin). Since the distance to the origin is monotonically
increasing between the point of minimal distance to the point of
maximal distance, the only such solution is the one described in
Theorem \ref{to self sim thm}, where the two half-lines are
parallel.\par
\noindent By Ben Andrews' results \cite{AndrewsClassification}, the
amount of curvature between $r_\text{min}$ and $r_\text{max}$ is at
least $\pi / 2$. Thus the only possibility (besides a network as in
Theorem \ref{to self sim thm}) to form such a loop is that on the
shorter arc we have only one minimum point and on the longer arc we
have two minima and one maximum. Let us parameterize this family of
solutions by $r_\text{min}$. Note that for a given $r_\text{min}$,
there are in general two corresponding values $\bar{\eta}\leq \eta_0$
and $\tilde{\eta} \geq \eta_0$.  To show existence and uniqueness we
only have to show that the amount of total curvature $K(r_\text{min})$
attains the value $\frac{4\pi}{3}$ only once. To do this write
$$K(r_\text{min})= 2\int_{r_\text{min}}^{r_\text{max}} \kappa\, ds + 4
\Psi(\eta) \equiv 2 \Theta(\rho) + 4 \Psi(\eta)\ ,$$ where
$\eta=\eta(r_\text{min})$ and $\rho=r_\text{max}/r_\text{min}$. Ben
Andrews has shown that $\lim_{\rho \rightarrow
  1}\Theta(\rho)=\frac{\pi}{\sqrt{2}},\ \lim_{\rho \rightarrow
  \infty}\Theta(\rho)=\frac{\pi}{2}$ and that $\Theta$ is monotone in
$\rho$, i.e. $\Theta(r_\text{min})$ is monotonically increasing in
$r_\text{min}$. Note also that in the interval $[1, \eta_0]$, $\bar{\eta}$
is monotonically increasing in $r_\text{min}$. From the proof of
existence and uniqueness of the completely symmetric lens, we know
that $\Psi(\eta)$ is strictly monotonically increasing in $\eta$ until
it attains the value $\pi/3$ and then it stays above this value. This
yields that
$$ \lim_{r_\text{min}\rightarrow 0}K(r_\text{min}) = 2 \lim_{\rho
  \rightarrow \infty} \Theta(\rho) = \pi\ . $$ Restricting to
$\bar{\eta}\leq\eta_0$, we see that $K(r_\text{min})$ is strictly
increasing until $\bar{\eta}=\eta_0$. For $\bar{\eta}=\eta_0$ we have
$K\geq \pi + 4\Psi(\eta_0)\geq \pi + 4\times 0.72 > 4\pi/3$. For
$\tilde{\eta}\geq \eta_0$, $\Psi$ is increasing or bounded below by
$\pi/3$ and we obtain the same lower bound. Thus $K$ attains the value
$4\pi/3$ exactly once. This implies existence and uniqueness in the
symmetric case.

\noindent {\bf Nonexistence in the Asymmetric Case:} Assume that an
asymmetric self-similarly shrinking generalized lens-shaped network
exists, i.\,e. a network which has no axis of reflection. Let us
suppose that the short arc is asymmetric, and thus also the long arc.
Once more we can assume that we have the same configuration of minima
and maxima as before. Let $\Sigma(\eta)$ be the total curvature of the
short asymmetric arc. Using the considerations in Section
\ref{selfsimilar sec} and Lemma \ref{bd for psieta0} we see that
$$\Sigma(\eta)\geq \Psi(\eta_0)>\frac{\pi}{6}\ . $$
Thus we get for the total curvature
$$  K(r_\text{min}) = 2 \Sigma(\eta) + 2\Theta(\rho) > \frac{\pi}{3} +
\pi = \frac{4\pi}{3} ,$$ which yields a contradiction.
\end{proof}
\end{appendix}

\bibliographystyle{amsplain}

\newpage
\end{document}